\documentclass[11pt,a4paper]{article}

\usepackage[T1]{fontenc}
\usepackage[utf8]{inputenc}
\usepackage{lmodern}
\usepackage{amsmath,amssymb,amsfonts}
\usepackage{amsthm}
\usepackage{mathrsfs}
\usepackage{graphicx}
\usepackage{xcolor}
\usepackage[hidelinks]{hyperref}
\usepackage{authblk} 
\usepackage{bm}      
\usepackage{enumerate}
\usepackage{geometry}
\theoremstyle{plain}
\newtheorem{theorem}{Theorem}[section]
\newtheorem{proposition}[theorem]{Proposition}

\theoremstyle{definition}
\newtheorem{definition}[theorem]{Definition}

\theoremstyle{remark}
\newtheorem{remark}{Remark}

\title{Non-local Boundary Value Problems, stochastic resetting and Brownian motions on graphs}

\author[1]{Stefano Bonaccorsi\thanks{stefano.bonaccorsi@unitn.it}}
\author[2]{Fausto Colantoni\thanks{fausto.colantoni@uniroma1.it (corresponding author)}}
\author[3]{Mirko D'Ovidio\thanks{mirko.dovidio@uniroma1.it}}
\author[4,5]{Gianni Pagnini\thanks{gpagnini@bcamath.org}}

\affil[1]{Department of Mathematics, University of Trento, Trento, Italy}
\affil[2]{Department of Basic and Applied Sciences for Engineering, Sapienza University of Rome, Rome, Italy}
\affil[3]{Department of Statistical Sciences, Sapienza University of Rome, Rome, Italy}
\affil[4]{BCAM - Basque Center for Applied Mathematics, Bilbao, Basque Country, Spain}
\affil[5]{Ikerbasque - Basque Foundation for Science, Bilbao, Basque Country, Spain}

\date{} 

\begin{document}
	
	\maketitle
	
	\begin{abstract}
		We consider dynamic boundary conditions involving non-local operators. Our  analysis includes a detailed description of such operators together with their relations with random times and random (additive) functionals. We provide some new characterizations for the boundary behaviour of the Brownian motion based on the interplay between non-local operators and boundary value problems. Our main focus is on Feller-Wentzell diffusions with jumps (resetting/restart). We first consider the instructive case of the real line, then we extend our results on star graphs with trapping points or repulsive vertices.
	\end{abstract}
	
	\section{Introduction}
	
	\subsection{Presentation of our results}
	In the present work we provide rigorous results for non-local boundary value problems (NLBVPs) on networks. We consider the problem to find a solution together with its probabilistic characterization to
	\begin{equation}
		\left\lbrace
		\begin{array}{ll}
			\displaystyle \dot{u}(t,\mathsf{x})= \Delta u(t,\mathsf{x}) & t>0, \; \mathsf{x} \in \mathcal{G}\setminus \{v\}, \\
			\displaystyle \eta \mathfrak{D}^\Psi_t u (t, \mathsf{x}) + \sum_{e \in \mathcal{E}} \rho_e\mathbf{D}^\Phi_{x-} u_e(t, \mathsf{x}) = 0 & t>0, \; \mathsf{x} \in \{v\},\\ 
			\displaystyle u(0,\mathsf{x})=f(\mathsf{x}) & \mathsf{x} \in \mathcal{G},
		\end{array}
		\right .
		\label{NLBVPgraphIntro}
	\end{equation}
	where $\mathcal{G}$ is a metric graph (star graph). The analysis on network can be carried out by considering union of star graphs.  The analysis on star graphs can be obtained from the half line via equivalence (see Theorem \ref{thm:equiv})). The non-local operators appearing above are given by the Caputo-D\v{z}rba\v{s}jan (type) derivative $\mathfrak{D}^\Psi_t$ and the Marchaud (type) derivative $\mathbf{D}^\Phi_{x-}$ where $\Psi$ and $\Phi$ are Bernstein symbols characterizing the operators. The probabilistic representation of the solution to \eqref{NLBVPgraphIntro} is given in Theorem \ref{thm:Xpallinodelay} after the presentation of the preparatory results. This representation on star graphs can be written in terms of a Brownian motion on $[0, \infty)$ reflected at zero for which we show that the spatial condition controls an additive part acting on the reflecting Brownian motion whereas the time condition introduces a time change acting on the local time (at zero) of the reflecting Brownian motion (see Theorem \ref{fracdyntm}). The additive part $\mathsf{L_{\gamma_t}}$ (see formula \eqref{process}) pushes away from zero the process and the time change $S_t^{-1}$ (see formula \eqref{TimeChangeS}) forces the process to stop for a random amount of time. Due to the right-continuity, the stop occurs immediately after the jump. We underline that the process stops for an independent amount of time (see Theorem \ref{thm:HTstar}) and jumps for an independent distance from the origin (the discussion around formula \eqref{process}). Such a behavior can be associated with a well-known dynamic often referred to as restart or resetting depending on the literature. Moreover we underline the following facts: i) if $\rho_e = 0$ for avery $e \in \mathcal{E}$, then $\Psi$ controls the holding times at the star vertex and $\Psi^\prime (0)=\infty$ means that such a vertex is trapping for the Brownian motion on $\mathcal{G}$ (see formula \eqref{meanHTgraph}); if $\eta=0$, then $\Phi$ controls the jumps away from the star vertex and the jumps can be associated with a subordinator with finite or infinite activity (local or non-local ${\bf D}^\Phi_{x-}$). We also underline that $\rho_e>0$ for some $e \in \mathcal{E}$ says that the star vertex is repulsive (we have jumps) whereas $\eta>0$ says that we may have finite holding times with infinite mean.  \\
	
	In Section 2 we introduce some basic facts and notations about subordinators and inverses. We discuss the relations between the operators we deal with and the associated processes. In particular, the left and right Marchaud (type) derivatives $\mathbf{D}^\Phi_{x-}$ and $\mathbf{D}^\Phi_{x+}$ can be respectively regarded as the generator and its adjoint for a subordinator with symbol $\Phi$. From this point of view, our analysis completely agrees with the boundary conditions introduced by Feller in \cite{feller1952parabolic}. Then, we discuss the heat equation on the half-line with non-local boundary conditions in Section 3. We provide the main result on the real line (Theorem \ref{fracdyntm}) which generalizes the work \cite{d2022fractionalsticky} for the time operator and the work \cite{feller1952parabolic} for the space operator. We also provide a discussion on the probabilistic representation of the solution by extending the results given in \cite{d2022fractional,d2022fractionalsticky} and the probabilistic representation given by It\^{o} and McKean in \cite{ito1996diffusion, ito1963brownian}. In Section 4 we generalize our results on metric graphs and provide our main result on graphs (Theorem \ref{thm:Xpallinodelay}) concerned with \eqref{NLBVPgraphIntro}. In Section 5 we discuss on interpretations, applications and extensions of our results. 
	
	\subsection{Discussion of related results and the existing literature}
	
	Let us consider the NLBVP on a metric space $E$,
	\begin{equation}
		\label{NLBVP}
		\left\lbrace
		\begin{array}{ll}
			\displaystyle \dot{u} = A u & t>0,\; x \in E \setminus \partial E,\\
			\displaystyle v= Tu, \; \mathfrak{D}^\Psi_t v = \dot{v} * \overline{\Pi}^\Psi & t>0, \; x \in \partial E,\\
			\displaystyle \mathfrak{D}^\Psi_t v + \Phi(B) u = 0 & t>0,\; x \in \partial E,\\
			\displaystyle u(0,x) = f(x) & x \in \overline{E},
		\end{array}
		\right .
	\end{equation}
	for a good initial datum $f$ and some operators $A,B$ for which $D(\Phi(B)) \subset D(B) \subset D(A)$. The symbol $\mathfrak{D}^\Psi_t$ denotes a convolution-type operator already introduced. Suppose that $E \subseteq \mathbb{R}^d$ is a smooth  domain and denote by $T$ the trace operator properly defined with respect to a finite Borel measure on $\partial E$. In general, the operator $\Phi(B)$ can be associated with a boundary motion. On $\mathbb{R}$, the problem is studied with a transmission condition at the origin, in order to characterize the non-local skew and non-local skew sticky Brownian motions \cite{colantoni2025non}. For a discussion in case $\Phi(B)=B$ is local, the reader can consult \cite{d2022fractionalsticky} in dimension $d=1$ and \cite{d2024fractionalsticky} in dimension $d>1$. Consider the case in which $(B, D(B))$ generates a Markov process on $\partial E$ and $\Phi(B)$ is a non-local operator defined via the Bernstein symbol $\Phi$ (Phillips' representation for instance). In case $\partial E$ is a compact manifold without boundary, the probabilistic reading of the non-local Cauchy problem (NLCP)
	\begin{equation}
		\label{NLCP}
		\left\lbrace
		\begin{array}{ll}
			\displaystyle \mathfrak{D}^\Psi_t u + \Phi(B) u = 0, & t>0,\; x \in \partial E,\\
			\displaystyle u(0,x) = f(x), & x \in \partial E,
		\end{array}
		\right .
	\end{equation}
	can be given in terms of a time-changed Markov process where the time change is obtained via subordinate semigroup associated with $\Phi(B)$ and an independent non Markov random time associated with $\mathfrak{D}^\Psi_t$. Roughly speaking the operator $\Phi(B)$ introduces jumps whereas the operator $\mathfrak{D}^\Psi_t$ introduces intervals of constancy. In turn, the non-local boundary condition in \eqref{NLBVP} introduces a behaviour near the boundary with a very tricky characterization.

	Depending on $d$, a process on $E$ may have interesting properties. For example a Brownian motion on $E \subset \mathbb{R}$ can hit its starting point infinitely many times. As $d>1$, there are many results on the visited points of a Brownian motion. For $d=2$ every point is almost surely not visited and despite of this, there exists a random set of points visited infinitely often. In both cases ($d=1$ and $d=2$) such sets are big and uncountable. For the process on  the metric space $E$ associated with \eqref{NLBVP} we focus on the set of boundary points with their hitting times. Although the problem \eqref{NLBVP} for $d>1$ is interesting, the case $E=[0, \infty)$ turns out to be very instructive. The boundary point $\{0\}$ has zero Lebesgue measure and represents a repulsive point as prescribed by $\Phi(B)$. On the other hand the operator $\mathfrak{D}^\Psi_t$ introduces random holding times on the boundary if $\Phi(B)=B$ or alternatively on a random point in the interior, after the jump, in case $\Phi(B)$ is non-local. \\

	In the present work we respectively focus on the half-line 
	\begin{align*}
		E=[0, \infty), \quad \partial E = \{0\}
	\end{align*}
	and the metric graph obtained from the collection of half-lines
	\begin{align*}
		E = \bigsqcup \, [0, \infty).
	\end{align*}
	We study the non-local dynamic boundary condition 
	\begin{align*}
		\mathfrak{D}^\Psi_t v + \Phi(B) u = 0
	\end{align*}
	where
	\begin{align*}
		\Phi(B) u = - \mathbf{D}^\Phi_{x-} u, \quad \textrm{on} \quad (0, \infty) \times \{0\}
	\end{align*}
	for the problem on the positive half-line and
	\begin{align*}
		\Phi(B) u = - \sum_{e \in \mathcal{E}}\, \rho_e \mathbf{D}^\Phi_{x-} u_e, \quad \textrm{on} \quad (0,\infty) \times \{v\} 
	\end{align*}
	for the problem on $\mathcal{G}$, the star graph.\\

	The intuition and the first rigorous results on boundary conditions associated with jumps away from the boundary have been given in \cite{feller1952parabolic, ito1963brownian} in case $E=[0, \infty)$. However, the connection with non-local operators has a recent development. The NLCP \eqref{NLCP} on compact domains has been investigated by many researcher in the last decades. In case of Caputo-D\v{z}rba\v{s}jan derivative in \cite{MNV09} the authors introduced the Fractional Cauchy Problem (FCP) on bounded domains and previously in \cite{MaLuPa2001} the space-time fractional diffusion equation on the real line has been investigated. The general problem has been studied in \cite{kochubei2011general} and \cite{toaldo2015convolution}. Concerning the NLBVP, in the literature there are clear references on dynamic boundary conditions. The non-local dynamic boundary condition as in \eqref{NLBVP} has been introduced, to the best of our knowledge, in \cite{d2024fractionalsticky} for bounded domains and in \cite{d2022fractional,d2022fractionalsticky} for the real line. We refer to the book \cite{bookMorPer} for the discussion on the visited points of the Brownian motion.\\

	\subsection{Motivations} Among the motivations driving our investigation into these models, lies a number of potential applications, from financial models to models of traffic for example. Specifically, dynamic non-local conditions prove invaluable for constructing a model of sticky expectations, wherein investors adjust their beliefs at a deliberately gradual pace. Meanwhile, spatial boundary conditions characterized by abrupt data jumps find relevance, for instance, in addressing structural breaks within financial systems. To further expand the scope of our research, we have extended our results to star graphs, offering a framework for networks with distinct behaviours. For example, traffic models with different holding times for given nodes. We recently obtained interesting applications also in \cite{ColDovTav} where NLBVPs on networks have been considered for the statistical analysis of earthquakes.\\

	\section{Preliminaries and notations} 
	\label{sec:PrelNot}		
	We recall some basic facts and introduce some notations. In order to streamline the notation as much as possible we write
	\begin{align*}
		\dot{u} = \frac{\partial u}{\partial t}, \quad u^\prime = \frac{\partial u}{\partial x}, \quad u^{\prime \prime} = \frac{\partial^2 u}{\partial x^2}.
	\end{align*}
	Moreover, we write
	\begin{align*}
		W^{1,p}(0, \infty) = \{ u \in L^p(0, \infty):\;  u^\prime \in L^p(0, \infty) \}
	\end{align*}
	and
	\begin{align*}
		W^{1,p}_0 (0, \infty) = \{ u \in W^{1,p}(0, \infty):\, u(0)=0\}
	\end{align*}
	with $p \in [1,\infty]$ for the Sobolev spaces. For the reader's convenience we list below some symbols:
	\begin{itemize}
		\item[-] $H^\Phi$, $H^\Psi$ are independent subordinators with symbols $\Phi$, $\Psi$;
		\item[-] $L^\Phi$, $L^\Psi$ are their inverse processes;
		\item[-] $\mathsf{L}$ is the remaining time in a location for the inverse of a subordinator/overshoot of the subordinator;
		\item[-] $\gamma(X)$ is the local time $\gamma$ of the process $X$. If no otherwise specified $\gamma=\gamma(B^+)$;
		\item[-] $B^+$ is a reflected Brownian motion;
		\item[-] $B^\bullet = B^+ + \mathsf{L} \circ \gamma$ is driven by (space) non-local boundary conditions;
		\item[-] $B^\star = B^\bullet \circ S^{-1}$ is driven by (space/time) non-local boundary conditions.
	\end{itemize}

	\subsection{Subordinators and symbols}
	Let $H^\Phi =\{H_t^\Phi , \ t \geq 0\}$ be a subordinator (see \cite{bertoin1999subordinators} for details). Then, $H^\Phi$ can be
	characterized by the Laplace exponent $\Phi$, that is,
	\begin{align}
		\label{LapH}
		\mathbf{E}_0[\exp(-\lambda H_t^\Phi)]=\exp(-t\Phi(\lambda)), \quad \lambda \geq 0.
	\end{align}
	We denote by $\mathbf{E}_x$ the expected value with respect to $\mathbf{P}_x$ where $x$ is the starting point. Since $\Phi(\lambda)$ is the Laplace exponent of a subordinator, it is uniquely characterized by the pair of non-negative real numbers $(\kappa,d)$  and by the L\'evy measure $\Pi^\Phi$ on $(0,\infty)$ such that $\int_0^\infty (1 \wedge z) \Pi^\Phi(dz) < \infty$.
	For the symbol $\Phi$, the following L\'evy-Khintchine representation holds (\cite[Theorem 3.2]{schilling2012bernstein})
	\begin{align}
		\label{LevKinFormula}
		\Phi(\lambda)=\kappa+d \lambda + \int_0^\infty (1-e^{-\lambda z}) \Pi^\Phi(dz), \quad \lambda>0
	\end{align}
	where the killing rate $\kappa$ and the drift coefficient $d$ are given by
	\begin{align*}
		\kappa=\lim_{\lambda \to 0} \Phi(\lambda), \quad d=\lim_{\lambda \to \infty} \frac{\Phi(\lambda)}{\lambda}.
	\end{align*}
	The symbol  $\Phi$ is a Bernstein function (non-negative, non-decreasing and continuous, see for example \cite{schilling2012bernstein}) uniquely associated with $H^\Phi$ (\cite[Theorem 5.2]{schilling2012bernstein}). For the reader's convenience we also recall that
	\begin{align}
		\label{tailSymb}
		\frac{\Phi(\lambda)}{\lambda}=d+ \int_0^\infty e^{-\lambda z} \overline{\Pi}^\Phi(z) dz, \quad \overline{\Pi}^\Phi(z)=\kappa+\Pi^\Phi(z,\infty)
	\end{align}
	where $\overline{\Pi}^\Phi$ is the so called \textit{tail of the L\'evy measure $\Pi^\Phi$}. We also define the process $L^\Phi=\{L_t^\Phi,\ t \geq 0\}$, with $L_0^\Phi=0$,  as the inverse of $H^\Phi$, that is
	\begin{align*}
		L_t^\Phi = \inf \{s > 0\,:\, H_s^\Phi >t \}, \quad t>0.
	\end{align*}
	
	{\sc Setting: } In this paper we only consider symbols for which $\kappa=0$ and $d=0$. \\
	
	Further on we will also introduce the processes $H^\Psi=\{H_t^\Psi, t \geq 0\}$, independent of $H^\Phi$, and the inverse $L^\Psi=\{L_t^\Psi, t \geq 0\}$ of $H^\Psi$ with symbol 
	\begin{align*}
		\Psi(\lambda):=\int_0^\infty (1-e^{-\lambda z}) \Pi^\Psi(dz), \quad \lambda >0.
	\end{align*}
	We always assume that $H^\Phi$ and $H^\Psi$ are independent subordinators.\\

	{ \sc Activity:} Let $\Pi^\Phi$ be the L\'{e}vy measure of the subordinator $H^\Phi$. We recall that $H^\Phi$ has finite (respectively infinite) activity if $\Pi^\Phi(0, \infty)< \infty$ (respectively $\Pi^\Phi(0, \infty) =  \infty$). The associated L\'{e}vy measure plays a role in the definition of non-singular (respectively singular) non-local operator we will deal with. \\

	In Section \ref{sec:resetting} we consider subordinators with finite activity and discuss the spacial case of stochastic resetting. Further on,  we mainly focus only on strictly increasing subordinators with infinite activity and zero drift (see \cite[Theorem 21.3]{sato1999levy}). Thus, the inverse process $L^\Phi$ turns out to be a continuous process. In particular, $H^\Phi$ may have jumps and the inverse $L^\Phi$ has non decreasing paths. Notice that, an inverse process can be regarded as an exit time for $H^\Phi$. By definition, we also have
	\begin{align}
		\label{relationHL}
		\mathbf{P}_0(H_t^\Phi < s) = \mathbf{P}_0(L_s^\Phi>t), \quad s,t>0.
	\end{align}
	
	{\sc On the involved functions I:}
	We deal here under the assumption that $H^\Phi$ admits density in order to provide a clear picture of the space of functions to be considered below. Notice that the density of $L^\Phi$ always exists under the setting previously introduced (see \cite[Theorem 3.1]{meerschaert2008triangular}). Let us introduce $h,l$ for which 		\begin{align*}
		\mathbf{P}_0(H_t^\Phi \in I) = \int_I h(t, x)dx,\quad  \mathbf{P}_0(L_t^\Phi \in I) = \int_I l(t, x)dx ,
	\end{align*}
	for a given set $I \subset (0,\infty)$. From \eqref{LapH}, we have that
	\begin{align}
		\label{Laph}
		\int_0^\infty e^{-\lambda x} h(t,x) dx= e^{-t\Phi(\lambda)}, \quad \lambda >0
	\end{align}
	and from \cite[formula (3.13)]{meerschaert2008triangular}, we get
	\begin{align}
		\label{LapL}
		\int_0^\infty e^{-\lambda t} l(t,x) dt=\frac{\Phi(\lambda)}{\lambda} e^{-x \Phi(\lambda)}, \quad \lambda >0.
	\end{align}
	By using the initial value theorem for the Laplace transform, we observe that, from \eqref{Laph},
	\begin{align}
		\label{limithzero}
		h(t,0)=\lim_{\lambda \to \infty} \lambda e^{-t \Phi(\lambda)}
	\end{align}
	only in case $\forall\, t$, $h(t, \cdot)$ is bounded. For example, $h(t,0)<\infty$. The limit above depends on $\Phi(\lambda)$ and the time variable $t$. This fact is known as \emph{time dependent property} (see Definition 23.1 in \cite{sato1999levy}).  Thus, the fact that $h(t,x)=0$ as $x \leq 0$ plays a special role in our discussion. Indeed, this condition has non trivial consequences and in general, for some $\Phi$ and some $t>0$,
	\begin{align*}
		h(t, \cdot) \notin W^{1,1}_0(0, \infty).
	\end{align*}
	The gamma subordinator with $\Phi(\lambda) = \ln (1+\lambda)$ is an example (see \cite{colantoni2021inverse}).\\

	\subsection{Marchaud (type) operators}
	\label{sec:Marchaud}
	For a continuous (causal) function $u$ on $\mathbb{R}$ extended with zero on the negative part of the real line, that is $u(x)=0$ if $x \leq 0$, we define the right Marchaud (type) derivative
	\begin{align}
		\label{Marchaudright}
		\mathbf{D}_{x+}^\Phi u(x)=\int_0^\infty (u(x) - u(x-y))\Pi^\Phi(dy)
	\end{align}
	and the left  Marchaud (type) derivative
	\begin{align}
		\label{Marchaudleft}
		\mathbf{D}_{x-}^\Phi u(x) = \int_0^\infty (u(x) - u(x+y))\Pi^\Phi(dy).
	\end{align}
	If $\Pi^\Phi$ is the L\'evy measure associated to a stable subordinator, formulas \eqref{Marchaudright} and \eqref{Marchaudleft} respectively coincide with the right and the left Marchaud derivatives, usually denoted by $\mathbf{D}_{x+}^\alpha$ and $\mathbf{D}_{x-}^\alpha$ respectively. The reader can consult the famous book \cite[formula (5.57) and (5.58)]{samko1993fractional} or the paper \cite[section 6.1]{ferrari2018weyl} for a recent discussion. The operators $\mathbf{D}_{x+}^\Phi$ and $\mathbf{D}_{x-}^\Phi$ can be defined on different spaces depending on $\Pi^\Phi$. For a general definition, given the symbol $\Phi$, we consider $u$  bounded and locally Lipschitz continuous, then
	\begin{align}
		\label{stimaMarchaud}
		\vert \mathbf{D}_{x\pm}^\Phi u(x) \vert &\leq \int_0^1 \vert u(x) -u(x\mp y) \vert \Pi^\Phi(dy) + \int_1^\infty \vert u(x) -u(x\mp y) \vert \Pi^\Phi(dy) \notag \\
		&\leq K \int_0^1 y \Pi^\Phi(dy) + 2 \vert\vert u \vert \vert_\infty \int_1^\infty \Pi^\Phi(dy) \notag \\
		&\leq(K+2 \vert \vert u \vert \vert_\infty) \int_0^\infty (1 \wedge y) \Pi^\Phi(dy) <\infty.
	\end{align}
	Indeed, the Lipschitz property for a positive constant $K>0$ holds in the first integral and the boundedness of $u$ holds in the second integral. Since $\int (1 \wedge z) \Pi^\Phi(dz) < \infty$, then the last inequality directly emerges.

	We provide the following two examples in case of explicit representations for $\Pi^\Phi$: 
	\begin{itemize}
		\item[-] The case of stable subordinator with $\Phi(\lambda)=\lambda^\alpha$ for $\alpha \in (0,1)$ and $\Pi^\alpha(dy)=\frac{\alpha}{\Gamma(1-\alpha)} \frac{dy}{y^{\alpha+1}}$. The operators \eqref{Marchaudright} and \eqref{Marchaudleft} are therefore defined for locally $\gamma-$H\"{o}lder continuous functions with $\gamma > \alpha$ (it can be easily checked by calculation);
		\item[-] The case of gamma subordinator with $\Phi(\lambda) = a\ln (1+\lambda/b)$ for $a,b>0$ and $\Pi^\Phi(dy)=a \frac{e^{-by}}{y} dy$. This case may be a bit demanding. This is due to the \textit{time dependent continuity} of $h(t,x)$ at $x=0$. The operator \eqref{Marchaudright} is defined for $\gamma-$H\"{o}lder continuous functions with $\gamma >0$ (see \cite{colantoni2021inverse} for details).
	\end{itemize}

	We also underline the following fact. The operator $\mathbf{D}_{x+}^\Phi$ can be obtained by the Phillips' representation (\cite{phillips1952generation}) in the set of functions extended with zero on the negative part of the real line. Indeed, for the shift semigroup $\mathcal{S}_y u(x) = u(x-y)$ we have the representation 
	\begin{align*}
		\mathbf{D}_{x+}^\Phi u(x)=\int_0^\infty (u(x) -\mathcal{S}_y u(x))\Pi^\Phi(dy)
	\end{align*}
	for which 
	\begin{align}
		\label{LapMarchaud}
		\int_0^\infty e^{-\lambda x} \mathbf{D}_{x+}^\Phi u(x) dx= \Phi(\lambda) \int_0^\infty e^{-\lambda x} u(x)\,dx
	\end{align}
	where we used \eqref{LevKinFormula} with $\kappa=0$ and $d=0$ in order to obtain $\Phi(\lambda)$. The representation given by Phillips well agrees with a Riemann-Liouville (type) operator described in Section \ref{Sec:Appendix} and denoted by $\mathcal{D}^\Phi_{(0,x)}$.  \\

	We notice that if $u \in W^{1,\infty}(0, \infty)$, then the Marchaud (type) derivatives \eqref{Marchaudright} and \eqref{Marchaudleft} are well defined almost everywhere. Indeed, with the first inequality of \eqref{stimaMarchaud} at hand, we use the fact that $u$ is essentially bounded and locally Lipschitz almost everywhere (its derivative is bounded).

	\subsection{Caputo-D\v{z}rba\v{s}jan (type) operators}
	
	The well-known Caputo-D\v{z}rba\v{s}jan derivative is associated with $\Phi(\lambda)=\lambda^\alpha$, $\alpha \in (0,1)$. It has been introduced independently in \cite{caputoBook, CapMai71,CapMai71b} and \cite{Dzh66, DzhNers68}. 
	
	Let $M>0$ and $w\geq 0$. Let $\mathcal{M}_w$ be the set of (piecewise) continuous function on $[0, \infty)$ of exponential order $w$ such that $\vert u(t)\vert \leq M e^{wt}$. Let $u \in \mathcal{M}_0 \cap C[0,\infty)$ with $u^\prime \in \mathcal{M}_0$. Then we define the Caputo-D\v{z}rba\v{s}jan (type) derivative as the convolution
	\begin{align}
		\label{Caputo}
		\mathfrak{D}_{x}^\Phi u(x):=\int_0^x u^\prime (y) \overline{\Pi}^\Phi(x-y) dy
	\end{align}
	whose Laplace transform writes
	\begin{align}
		\label{LapCaputo}
		\int_0^\infty e^{-\lambda x}  \mathfrak{D}_{x}^\Phi u(x) dx= \Phi(\lambda) \widetilde{u}(\lambda)- \frac{\Phi(\lambda)}{\lambda} u(0) = \frac{\Phi(\lambda)}{\lambda} \big( \lambda \widetilde{u}(\lambda) - u(0) \big), \quad \lambda > 0
	\end{align}
	($\widetilde{u}$ denotes the Laplace transform of $u$). For a detailed discussion on the operator $\mathfrak{D}^\Phi_{x}$ see for example the recent works \cite{kochubei2011general, toaldo2015convolution, chen2017time}.
	
	\begin{proposition}
		Let the previous setting prevails. For the convolution \eqref{Caputo}, we have that 
		\begin{align}
			\label{normCaputo}
			\vert \vert \mathfrak{D}_{x}^\Phi u \vert \vert_p \leq \vert\vert u^\prime \vert \vert_p \left(\lim_{\lambda \to 0} \frac{\Phi(\lambda)}{\lambda}\right), \quad p \geq 1.
		\end{align}
	\end{proposition}
	\begin{proof}
		First we observe that the Laplace transform can be considered in order to obtain 
		\begin{align*}
			\lim_{\lambda \to 0} \int_0^\infty e^{-\lambda t} |\mathfrak{D}_{x}^\Phi u|^p\, dt = \| \mathfrak{D}_{x}^\Phi u \|_p^p.
		\end{align*}
		However, from the Young's inequality, we have
		\begin{align*}
			\| \mathfrak{D}_{x}^\Phi u \|_p \leq \| u^\prime \|_p \left( \int_0^\infty \overline{\Pi}^\Phi (y) dy \right)
		\end{align*}
		where only the last integral can be obtained via Laplace transform. Indeed,
		\begin{align*}
			\int_0^\infty \overline{\Pi}^\Phi (y) dy 
			= & \lim_{\lambda \to 0} \int_0^\infty e^{-\lambda y} \overline{\Pi}^\Phi (y) dy\\ 
			= & [\textrm{from }\eqref{tailSymb}\textrm{ with }\kappa=d=0 ]\\ 
			= & \lim_{\lambda \to 0} \frac{\Phi(\lambda)}{\lambda}, \quad \lambda>0.
		\end{align*}
	\end{proof}
	
	We remark that $\lim_{\lambda \to 0} \Phi(\lambda)/\lambda$ is finite only in some cases (see \cite{capitanelli2019fractional}). If $\lim_{\lambda \to 0} \Phi(\lambda)/\lambda$ is finite, then \eqref{normCaputo} gives a clear information about the cases of measurable or essentially bounded functions. In particular, we have that for $u^\prime \in L^1(0,\infty)$, if 
	\begin{align*}
		\lim_{\lambda \to 0} \frac{\Phi(\lambda)}{\lambda} < \infty,
	\end{align*}
	we can apply \eqref{normCaputo} and the derivative exists almost everywhere. By considering the well-known relation 
	\begin{align*}
		\mathfrak{D}_x^\Phi u(x)= \mathcal{D}_{(0,x)}^\Phi \big(u(x) - u(0) \big)
	\end{align*}
	(the definition of $\mathcal{D}_{(0,x)}^\Phi$ has been postponed in Section \ref{Sec:Appendix}) we have equivalence between Caputo and Riemann-Liouville derivatives as $u(0)=0$. We can check that, in this case, \eqref{LapCaputo} coincides with \eqref{LapMarchaud}. This entails a connection with causal functions, that is $\mathfrak{D}^\Phi_x u$ is well-defined for $u \in W^{1, \infty}(0, \infty)$.\\

	{\sc On the involved functions II:} Let us consider the class of functions for which
	\begin{align}
		\label{condMD}
		\exists\, M>0\,:\, | u^\prime(y) | \leq  M\, \frac{\Lambda^\Phi(dy)}{dy}
	\end{align} 
	where 
	\begin{align*}
		\Lambda^\Phi(dy) = \int_0^\infty \mathbf{P}_0(H_t^\Phi \in dy) dt
	\end{align*}
	is the potential measure for the subordinator $H^\Phi$ with symbol $\Phi$. Since $\Lambda^\Phi$ and $\overline{\Pi}^\Phi$ are associated Sonine kernels for which 
	\begin{align*}
		\int_0^x \overline{\Pi}^\Phi(x-y) \Lambda^\Phi(dy) =1
	\end{align*} 
	and
	\begin{align}
		| \mathfrak{D}^\Phi_x u(x)| \leq M \int_0^x \overline{\Pi}^\Phi(x-y) \Lambda^\Phi(dy),
	\end{align} 
	we get that $|\mathfrak{D}^\Phi_x u(x)|$ is uniformly bounded on $(0, \infty)$.\\

	{\sc On the symbol $\Phi$:} From \eqref{LapH} we can immediately see that $\mathbf{E}[H^\Phi_t] = t\, \Phi^\prime(0)$ where
	\begin{align}
		\Phi^\prime(0) = \lim_{\lambda \to 0} \frac{\Phi(\lambda)}{\lambda}.
		\label{PHIderivative}
	\end{align}
	Thus, formula \eqref{normCaputo} provides a link in the space $L^p$ between the potential \eqref{LapCaputo} and the moments of $H^\Phi$.

	\section{Motions on half-lines} 
	\label{sec:BC}
	
	\subsection{Non-local Boundary Value Problems}
	
	In this section we discuss the connection between subordinators and non-local boundary conditions. According with \cite{MeSt14} we introduce the regenerative set
	\begin{align*}
		\mathbf{M} = \{(t,\omega) \subset \mathbb{R} \times \Omega\,:\, t = H^\Phi_z(\omega) \textrm{ for some } z\geq 0\}
	\end{align*}
	where $H^\Phi_t=H^\Phi_t(\omega)$ is a pure jump subordinator with symbol $\Phi$. This definition does not agree with the definition given for example in \cite{Mais1983} and others but it allows a nice formulation of our problem. Moreover, we obtain a random closed set such that the right-hand portion of $\mathbf{M}$ is independent (and identically distributed) from the left-hand portion of $\mathbf{M}$ where the portions can be regarded as determined by a stopping time.  The random set $\mathbf{M}$ is a set of renewal points for the inverse $L^\Phi_t$. In particular, for a given path of the subordinator, the set $\mathcal{M}(\omega) = \{t \,:\, (t,\omega) \in \mathbf{M} \}$ is a countable union of the intervals $[H^\Phi_{z-}, H^\Phi_z)$.  For any $t\geq 0$ we define the last time of regeneration before $t$,
	\begin{align*}
		\mathsf{r}_t := \sup\{s \leq t\,:\, s\in \mathbf{M}(\omega) \}
	\end{align*}
	and the next time of regeneration after $t$, 
	\begin{align*}
		\mathsf{R}_t := \inf\{s > t\,:\, s \in \mathbf{M}(\omega) \}
	\end{align*}
	for which $\mathsf{r}_t \leq t \leq \mathsf{R}_t$. As noted in \cite{MeSt14} and in \cite[Section 2]{bertoin1999subordinators} we have that a.s. 
	\begin{align*}
		\mathsf{r}_{t-} = H^\Phi_{L^\Phi_t -} := H^{\Phi,-} \circ (L^\Phi_t) \quad \textrm{and} \quad \mathsf{R}_t = H^\Phi_{L^\Phi_t} := H^\Phi \circ L^\Phi_t,
	\end{align*}
	\(H^{\Phi,-}_t:=H^\Phi_{t-} \). We also introduce the age process 
	\begin{align*}
		\mathsf{A}_t := t - \mathsf{r}_t
	\end{align*}
	and the remaining lifetime 
	\begin{align*}
		\mathsf{L}_t := \mathsf{R}_t - t.
	\end{align*}
	%
	For a given $t$, the total amount of time the particle stops in its location is given by $\mathsf{A}_t + \mathsf{L}_t$. A particle that is already in its location according with $\mathsf{A}_t$ will move again according with $\mathsf{L}_t$.

	We now introduce the main object we deal with, that is the process \text{$B^\bullet = \{B^\bullet_t, t\geq 0\}$} written as
	\begin{align}
		\label{process}
		B^\bullet_t = B^+_t + \mathsf{L}_{\gamma_t}
	\end{align}
	where $\mathsf{L}_{\gamma_t}:= \mathsf{L} \circ \gamma_t$ and $\gamma_t=\gamma_t(B^+)$ is the local time at zero of the reflecting Brownian motion $B^+ = \{B^+_t, \, t\geq 0\}$ on $[0, \infty)$. The process $B^\bullet$ has been first introduced in \cite{ito1963brownian} in order to provide a probabilistic representation for a class of problems involving general boundary conditions. In this regard, a current reading can be given in terms of Marchaud (type) operators. The probabilistic reading of \eqref{process} says that $B^{\bullet}$ is a reflecting Brownian motion jumping away from the origin, with a length of the jump given by the subordinator $H^\Phi$ running with the clock $L^\Phi_{\gamma_t}$. Since $\Pi^\Phi(0,\infty)=\infty$, then in each time interval the number of jumps of $H^\Phi$ are infinite. A  detailed analysis of the sample paths has be given in \cite[section 12]{ito1963brownian}. We provide a discussion below.

	We remark that $L^\Phi \circ H^\Phi_t = t$ almost surely, whereas the study of composition $H^\Phi \circ L^\Phi_t$ is a bit demanding. We recall that (\cite[Proposition 2, section III.2]{bertoin1996levy}) 
	\begin{align*}
		\mathbf{P}_0( \mathsf{r}_{t-} < t = \mathsf{R}_t ) = 0.
	\end{align*}
	If $d = 0$ in \eqref{LevKinFormula}, then $\mathbf{P}_0(\mathsf{R}_t > t) = 1$ for every $t > 0$ (\cite[Theorem 4, section III.2]{bertoin1996levy}). Thus, for a zero drift subordinator, the passage above any given level is almost surely realized by a jump. Notice that $\mathsf{L}_{\gamma_t}$ is constant as $B^+ >0$.\\

	Now we focus on the problem with non-local dynamic boundary conditions and provide the probabilistic representation of the solution. With the discussions about the involved functions of Section \ref{sec:PrelNot} in mind, we introduce $D_0=D_1 \cap D_2$ where 
	%
	
	\begin{align*}
		D_1 = \left\lbrace \varphi : \forall\, x,\; \varphi(\cdot,x) \in W^{1,\infty}(0,\infty) \;\; \textrm{and} \;  \lim_{x \downarrow 0} \mathfrak{D}^\Psi_t \varphi(t,x) \; \textrm{exists} \right\rbrace,
	\end{align*}
	
	\begin{align*}
		D_2 = \left\lbrace \varphi : \forall\, t,\; \varphi(t,\cdot) \in W^{1,\infty}(0,\infty) \; \textrm{and} \;\; \lim_{x\downarrow 0} \mathbf{D}^\Phi_{x-} \varphi(t,x) \; \textrm{exists} \right\rbrace.
	\end{align*}
	
	Let us define the random time
	\begin{align}
		S_t = t + H^\Psi \circ (\eta L^\Phi_{\gamma_t}), \quad L^\Phi_{\gamma_t} := L^\Phi \circ \gamma_t
		\label{TimeChangeS}
	\end{align}
	and its inverse
	\begin{align*}
		S_t^{-1} = \inf\{ s>0\,:\, S_s>t\}.
	\end{align*}
	Before the main theorem on $[0,+\infty)$, we first explain what kind of integration we are considering.
	\begin{remark}
		Let $a(t)$ be a function from $[0,a]$ to $[0,\infty)$ which is non-decreasing and right-continuous, and satisfies $a(0)=0$, $a(\infty)=\lim_{t \to \infty} a(t)$.
		From \cite[Lemma 2.2, section V]{blumenthal}, for every nonnegative Borel measurable function $f$ on the positive half-line vanishing at infinity,
		\begin{align*}
			\int_{(0,\infty)} f(t) d a(t)= \int_0^\infty f(a^{-1} (t)) dt,
		\end{align*}
		where $a^{-1} (t)$ is the right inverse of $a(t)$. In this sense we interpret the integrals of Theorem \ref{fracdyntm}.
	\end{remark}
	
	We recall that, in this paper, \(\Phi\) and \(\Psi\) are two Bernstein functions with no drift. They are associated with two independent subordinators, \(H^\Phi\) and \(H^\Psi\), each with infinite activity.
	\begin{theorem}
		\label{fracdyntm}
		The solution $\upsilon \in C([0,\infty),(0,\infty)) \cap D_0$ to the problem
		\begin{align*}
			\begin{cases}
				\displaystyle 	\dot{\upsilon}(t,x)= \upsilon^{\prime \prime}(t,x) & t>0,\; x \in (0,\infty),\\
				\displaystyle 	\eta\, \mathfrak{D}_t^\Psi \upsilon(t,x) = -\mathbf{D}_{x-}^\Phi \upsilon(t,x), \; \eta \geq 0    & t>0,\; x=0,\\
				\displaystyle 	\upsilon(0,x)=f(x), \; f \in C[0,\infty) \cap L^\infty(0, \infty)  &x\in [0, \infty)
			\end{cases}
		\end{align*}
		has the probabilistic representation
		\begin{align}
			\label{NLstickyBpallino}
			\upsilon(t,x)=\mathbf{E}_x\left[f(B^\bullet \circ S_t^{-1})\right].
		\end{align}
	\end{theorem}
	We postpone the proof in Section \ref{proof1}.\\

	The lifetime of the process $B^\bullet \circ S^{-1}$ (which basically moves on a given path of a reflecting Brownian motion) is infinite, that is $R_0 \mathbf{1}=\infty$ and $\forall \, x$, $\upsilon(\cdot, x) \notin L^1(0,\infty)$ for a bounded initial datum $f$. However, there is no hope to find $\upsilon(\cdot,x) \in L^1(0,\infty)$ by also asking for $f \in L^1(0,\infty)$ for which we only have 
	\begin{align*}
		R_0^D f(x)= 2 \int_0^\infty (y\wedge x) f(y) dy < \infty 
	\end{align*}
	in the potential $R_\lambda f =R_\lambda^D f+ \bar{R}_\lambda f$.
	The last formula can be obtained from \eqref{Rtemp}. Then we ask for $f \in L^\infty$ in order to obtain $\upsilon (\cdot,x) \in L^\infty(0,\infty)$, $\forall x$.

	The previous result is concerned with the following particular cases:
	\begin{itemize}
		\item[$\eta=0$)] Non-local space conditions. The problem takes the form
		\begin{align*}
			\begin{cases}
				\displaystyle 	\dot{\upsilon}(t,x)= \upsilon^{\prime \prime}(t,x) & t>0,\; x \in (0,\infty),\\
				\displaystyle  \mathbf{D}_{x-}^\Phi \upsilon(t,x)=0   & t>0,\; x=0,\\
				\displaystyle 	\upsilon(0,x)=f(x), \; f \in C[0,\infty) \cap L^\infty(0, \infty)  &x\in [0, \infty).
			\end{cases}
		\end{align*}
		From the probabilistic point view, $S^{-1}_t$ becomes the inverse to $S_t = t$ and therefore, the solution $\upsilon$ has the representation
		\begin{align*}
			\upsilon(t,x) = \mathbf{E}_x \left[f(B^\bullet_t) \right].
		\end{align*}
		We observe that
		\begin{align}
			\label{integralMarchaud}
			-\mathbf{D}_{x-}^\Phi \upsilon(t,x) \big\vert_{x=0} 
			= & \lim_{x \to 0} \int_0^\infty \big(\upsilon(t,x+y) - \upsilon(t,x)\big) \Pi^\Phi(dy) \notag \\
			= & \int_0^\infty \big(\upsilon(t,y) - \upsilon(t,0)\big) \Pi^\Phi(dy)
		\end{align}
		which is the boundary condition introduced by W. Feller in  \cite{feller1952parabolic} and studied by K. It{\^o} and P. McKean  in \cite[section 12]{ito1963brownian}.\\
		
		We now discuss the special case of the $\alpha-$stable subordinator $H^\alpha$, characterized by the symbol
		\begin{align*}
			\Phi(\lambda)=\lambda^\alpha=\frac{\alpha}{\Gamma(1-\alpha)} \int_0^\infty (1-e^{-\lambda y}) \frac{dy}{y^{\alpha+1}}, \quad \alpha \in (0,1).
		\end{align*}
		Then, we write $\mathbf{D}_{x-}^\alpha$ in place of $\mathbf{D}_{x-}^{\Phi}$ which coincides with the well known Marchaud left derivative. Let us consider the following properties for the Marchaud derivatives
		\begin{align}
			\label{Marchaudalphazero+}
			\mathbf{D}_{x+}^\alpha u(x) &\to u(x) \quad \ \ \alpha \downarrow 0,\\
			\label{Marchaudalphaone+}
			\mathbf{D}_{x+}^\alpha u(x) &\to u^\prime(x) \quad \  \alpha \uparrow 1.
		\end{align}
		Such continuity properties hold for the Riemann-Liouville derivatives (see for example \cite[(2.2.5)]{kilbas2006theory}), we can adapt such results and obtain \eqref{Marchaudalphazero+} and \eqref{Marchaudalphaone+} just considering $u$ continuously differentiable with $u^\prime(x)$ that vanishes at infinity as $\vert x \vert^{\alpha-1-\varepsilon}$, $\varepsilon>0$ (see \cite[Section 5.4]{samko1993fractional}). The following analogue formulas hold true
		\begin{align}
			\label{Marchaudalphazero}
			\mathbf{D}_{x-}^\alpha u(x) &\to u(x) \quad \ \ \alpha \downarrow 0,\\
			\label{Marchaudalphaone}
			\mathbf{D}_{x-}^\alpha u(x) &\to -u^\prime(x) \quad \alpha \uparrow 1 .
		\end{align}
		Indeed, formula \eqref{MarchaudAdjoint} below which still holds for $\mathbf{D}_{x\pm}^\Phi u \in L^\infty(0,\infty)$, together with \eqref{Marchaudalphazero+} and \eqref{Marchaudalphaone+} say that \eqref{Marchaudalphazero} and \eqref{Marchaudalphaone} hold in some sense. Thus, roughly speaking,  as $\alpha \downarrow 0$ we get
		
		\begin{align}
			\label{killedBM}
			\begin{cases}
				\displaystyle 	\dot{\upsilon}(t,x)= \upsilon^{\prime \prime}(t,x) & t>0,\; x \in (0,\infty),\\
				\displaystyle  \upsilon(t,x)=0   & t>0,\; x=0,\\
				\displaystyle 	\upsilon(0,x)=f(x), \; f \in C(0,\infty) \cap L^\infty(0, \infty)  &x\in (0, \infty).
			\end{cases}
		\end{align}
		and for $\alpha \uparrow 1$, we get
		
		\begin{align}
			\label{reflectedBM}
			\begin{cases}
				\displaystyle 	\dot{\upsilon}(t,x)= \upsilon^{\prime \prime}(t,x) & t>0,\; x \in (0,\infty),\\
				\displaystyle  \upsilon^\prime (t,x)=0    & t>0,\; x=0,\\
				\displaystyle 	\upsilon(0,x)=f(x), \; f \in C[0,\infty) \cap L^\infty(0, \infty)  &x\in [0, \infty).
			\end{cases}
		\end{align}

		It is well-known that almost surely, for $\alpha=0$ the subordinator dies immediately whereas for $\alpha=1$, it is the elementary subordinator $t$ (\cite[Section 3.1.1]{bertoin1999subordinators}). The process \eqref{process} for $\alpha \downarrow 0$ becomes a killed Brownian motion as expected for the solution of \eqref{killedBM}. For $\alpha \uparrow 1$, the subordinator becomes $H_t^\alpha = t$ and the process \eqref{process} is a reflected Brownian motion that does not jump at the boundary point. It reflects and never dies. Indeed, it coincides with a reflected Brownian motion on $[0,\infty)$ as expected for the solution of \eqref{reflectedBM}.
		%

		\item[$\Phi=Id$)] Non-local dynamic conditions. For the problem
		\begin{align*}
			\begin{cases}
				\displaystyle 	\dot{\upsilon}(t,x)= \upsilon^{\prime \prime}(t,x) & t>0,\; x \in (0,\infty),\\
				\displaystyle 	\eta\, \mathfrak{D}_t^\Psi \upsilon(t,x) = \upsilon^\prime(t,x), \; \eta \geq 0   & t>0,\; x=0,\\
				\displaystyle 	\upsilon(0,x)=f(x), \; f \in C[0,\infty) \cap L^\infty(0, \infty)  &x\in [0, \infty)
			\end{cases}
		\end{align*}
		the solution $\upsilon$ takes the following probabilistic representation
		\begin{align}
			\upsilon(t,x)=\mathbf{E}_x\left[f(B^+ \circ (t+H^\Psi \circ \eta{\gamma_t})^{-1})\right].
			\label{Vbar}
		\end{align}
		This results well agrees with \cite[Theorem 3.1]{d2022fractionalsticky}. Here we have $H^\Phi_t=t$ a.s. implies that $\mathsf{L}_t = 0$ a.s. and $B^\bullet = B^+$. Moreover, $S_t$ equals in law $t + H^\Phi \circ \eta \gamma_t$. 
		
		\item[$\Psi=Id$)] Dynamic condition and non-local space conditions. We have
		\begin{align*}
			\begin{cases}
				\displaystyle 	\dot{\upsilon}(t,x)= \upsilon^{\prime \prime}(t,x) & t>0,\; x \in (0,\infty),\\
				\displaystyle 	\eta\, \dot{\upsilon}(t,x) = -\mathbf{D}_{x-}^\Phi \upsilon(t,x), \; \eta \geq 0   & t>0,\; x=0,\\
				\displaystyle 	\upsilon(0,x)=f(x), \; f \in C[0,\infty) \cap L^\infty(0, \infty)  &x\in [0, \infty).
			\end{cases}
		\end{align*}
		The boundary condition can be written as $\Delta \upsilon(t,x)\big\vert_{x=0}= -\mathbf{D}_{x-}^\Phi \upsilon(t,x) \big\vert_{x=0}$ by means of which we can underline the sticky nature of the motion.
		Indeed, $\mathfrak{D}_t^\Psi \upsilon= \partial_t \upsilon$ in Theorem \ref{fracdyntm} and $\upsilon$ satisfies the heat equation on $[0,\infty)$. The solution has the probabilistic representation
		\begin{align}
			\label{stickyBpallino}
			\upsilon(t,x)=\mathbf{E}_x\left[f(B^\bullet \circ (t+\eta L^\Phi_{\gamma_t})^{-1})\right]
		\end{align}
		where $B^\bullet_t$ is defined in \eqref{process} and $t+\eta L^\Phi_{\gamma_t} \stackrel{d}{=} t + \eta \gamma^\bullet_t$  must be considered in place of $S_t$ where $\gamma^\bullet_t = \gamma_t(B^\bullet)$. Indeed, $\gamma^\bullet_t \stackrel{d}{=} L^\Phi_{\gamma_t}$ with $\gamma_t= \gamma_t (B^+)$ (see \cite[Section 14]{ito1963brownian}) and $H^\Psi_t = t$ almost surely in case $\Psi=Id$. 
		
	\end{itemize}

		\begin{remark}
			({\it Non-local reflection on bounded domains}) We highlight the main issue in dealing with \( B^\bullet \circ S^{-1} \) in a bounded domain. Since the process exhibits jumps as it approaches the boundary, it is necessary to control the jumping measure to ensure that the process cannot jump outside the domain. A preliminary result for \(\alpha\)-stable L\'{e}vy processes can be found in \cite{Bogdan}. The authors are currently working to provide an explicit explanation through non-local operators at the boundary. On the other hand, the case of compactly supported measures (for the jumps), which leads to non-local non-singular operators, appears to be more manageable. For \( \Phi = \text{Id} \), a clear understanding of the problem is provided in \cite{d2022fractionalsticky} for smooth domains \( \Omega \subset \mathbb{R}^d \).
		\end{remark}

		\subsection{Sample paths description}
\begin{figure}[h!]
			\centering
			\includegraphics[width=10cm]{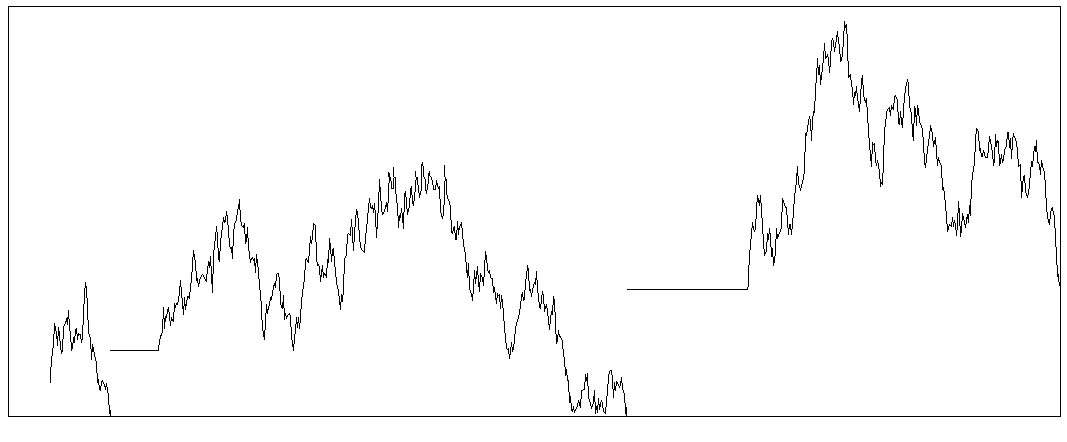} 
			\caption{A possible path for $B^\bullet \circ S^{-1}$. The case $\eta>0$. Since $\Psi \ne Id$ the random plateaus are given by $H^\Psi$ which is independent from $B^\bullet$.}
			\label{StickyBpallino}
		\end{figure}

		We recap briefly the dynamics of the processes introduced in Section \ref{sec:BC}. The process $B^\bullet$ defined in \eqref{process} is a reflected Brownian motion that jumps from $0$ inside $(0,\infty)$ in a random point given by the jumps of $H^\Phi$. This process turns out to be right-continuous since $H^\Phi \circ L^\Phi$ is the composition of the right-continuous subordinator $H^\Phi$ with its inverse (and continuous process) $L^\Phi$.
		\begin{figure}[h]
			\centering
			\includegraphics[width=10cm]{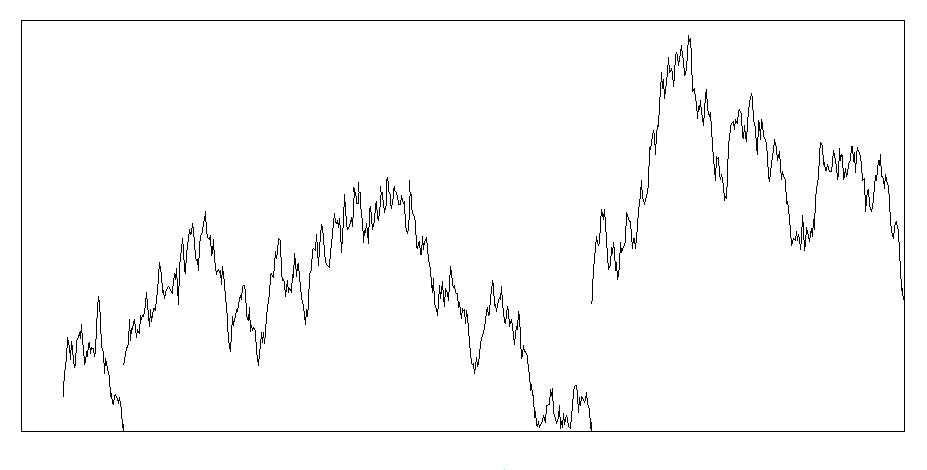} 
			\caption{A possible path for $B^\bullet$. The case $\eta=0$. The process is pushed away from the boundary point, it never hits $x=0$ and it jumps randomly in $(0,\infty)$ according with $\mathsf{L}_{\gamma_t}$.}
			\label{Bpallino}
		\end{figure}
		\begin{figure}[h]
			\centering
			\includegraphics[width=10cm]{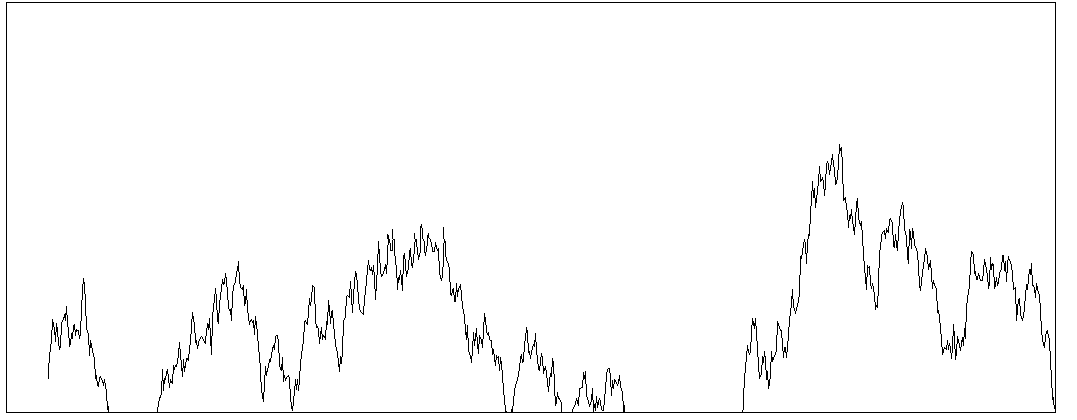} 
			\caption{A possible path for $B^+\circ \bar{V}^{-1}$. The case $\eta>0$ and $\Phi=Id$.}
			\label{onlysticky}
		\end{figure}
		The Sticky behavior is controlled by the non-local dynamic condition ($\eta >0$) and realized via time-change. From the construction we see that the reference path is given by the reflected Brownian motion. As $B^+$ hits the origin, then $\gamma_t$ increases in a such a way that the time-change runs slowly and slows down the process. The boundary condition $\eta \, \dot{\upsilon}(t,0)= -\mathbf{D}_{x-}^\Phi \upsilon(t,x) \big\vert_{x=0}$ leads to \eqref{stickyBpallino} which is the case of a (reflected) sticky process. Due to right-continuity of $H^\Phi$, determining an independent jump for $B^\bullet$, the time change $(t+\eta L^\Phi_{\gamma_t})^{-1}$ slows down $B^\bullet$ immediately after the jump. Then, an independent (exponential) holding time determines the sticky behavior at that state. The same arguments apply for the boundary condition $\eta\mathfrak{D}_t^\Psi \upsilon(t,x)\big\vert_{x=0}=-\mathbf{D}_{x-}^\Phi \upsilon(t,x) \big\vert_{x=0}\quad$ leading to \eqref{NLstickyBpallino}. The holding time after the jump is no more exponential, it may have heavy tailed distribution. Figure \ref{StickyBpallino} shows and example of these paths, a plateau occurs after a jump. We only underline once again that the plateaus are given by $H^\Psi$ whereas the jumps are given by $H^\Phi$ and they are both independent from $B^+$ which is the base process in the characterization we have proposed. Let us go into more detail to understand the paths of \(B^\bullet \circ S^{-1}\). Since \(L^\Phi \circ \gamma\) is a continuous process and \(H^\Psi\) is a pure jump process, we have \(\forall s>0\)
		\begin{align*}
			\Delta (H^\Psi \circ \eta L^\Phi \circ \gamma)_s:= (H^\Psi \circ \eta L^\Phi \circ \gamma)_s - (H^\Psi \circ \eta L^\Phi \circ \gamma)_{s-}= (\Delta H^\Psi) \circ (\eta L^\Phi \circ \gamma)_s.
		\end{align*}
		Then, \(\Delta (H^\Psi \circ \eta L^\Phi \circ \gamma)>0\) if and only if \(\eta L^\Phi \circ \gamma\) is in a jumping time of \(H^\Psi\). On the other hand, because \(S\), as \(H^\Psi\), is a right-continuous non decreasing process, we see
		\[t \in (S_{s-}, S_s] \iff S^{-1}_t=s,\]
		i.e., corresponding to the jumps of \(S\), we have intervals of constancy of \(S^{-1}\). Let us introduce the jumping times of \(S\)
		\[\mathfrak{J}_S:=\{s \geq 0 : \Delta S_s>0\}=\{ s \geq 0 : \,(\Delta H^\Psi) \circ (\eta L^\Phi \circ \gamma)_s>0\},\]
		where \(\Delta S_t:=S_t - S_{t-}\), and for \(s \in \mathfrak{J}_S\), we define \(a_s:=S_{s-}\) and \(b_s:=S_s=a_s+(\Delta H^\Psi) \circ (\eta L^\Phi \circ \gamma)_s\).  For \(t \notin \mathfrak{J}_S\), we have that \(B^\bullet \circ S^{-1}\) behave like \(B^\bullet\), since \(S\) is not jumping.
		\\
		For \(t \in \mathfrak{J}_S\), i.e. \(t \in (a_s, b_s]\), we have observed that \(S^{-1}_t=s\). So, we obtain
		\[B^\bullet_{S^{-1}_t} = B^\bullet_s \quad \forall t \in (a_s, b_s], \]
		which means that the process is constant on  \((a_s, b_s]\). From the right-continuity of \(B^\bullet\), we get
		\[B^\bullet_{S^{-1}_{a_s-}} = B^\bullet_{s-}=0, \quad B^\bullet_{S^{-1}_{a_s}} = B^\bullet_{s}=\xi_s>0,\]
		where \(\xi_s\) is the jump of \(B^\bullet\). We conclude that
		\[B^\bullet_{S^{-1}_t}= B^\bullet_s =\xi_s \quad \forall t \in (a_s, b_s], \]
		which guarantees that the process remains constant in the jump position.
		
		\subsection{Holding and renewal times}
		\label{sec:HTandRT}
		Let us consider the construction given in \cite{d2022fractionalsticky} and \cite{d2024fractionalsticky}. We respectively denote by $\{e_i,\, i \in \mathbb{N}\}$ and $\{\bar{e}_i := H^\Psi_{e_i},\, i \in \mathbb{N}\}$ the sequences of holding times (at the boundary) for the reflecting sticky Brownian motion, say $X$, and the partially delayed reflecting Brownian motion, say $\bar{X}$. We say that $\bar{X}$ is partially delayed in the sense that a delay effect only occurs on the boundary (the process behaves like a Brownian motion away from the boundary point $x=0$). Indeed,  $\bar{X}$ equals in law
		\begin{align*}
			B^+ \circ (t + H^\Psi \circ \eta \gamma_t)^{-1} \quad \textrm{where} \quad \gamma_t = \gamma_t(B^+), \quad \eta \geq 0.
		\end{align*} 
		The process $X$ corresponds to the case  $\Psi=Id$. In particular, $X$ is a reflected sticky Brownian motion driven by the problem
		\begin{align*}
			\begin{cases}
				\displaystyle \dot{\upsilon}(t,x)= \upsilon^{\prime \prime}(t,x) & t>0,\, x \in (0,\infty), \\
				\displaystyle \eta\, \dot{\upsilon}(t,x) = \upsilon^\prime(t,x) & t>0,\, x=0,\\
				\displaystyle \upsilon(0,x)=f(x), \; f \in C(0, \infty) \cap L^\infty(0, \infty) & x \in [0, \infty),
			\end{cases}
		\end{align*}
		for $\eta \geq 0$. That is,
		\begin{align*}
			\upsilon(t,x) = \mathbf{E}_x\left[f(X_t)\right] = \mathbf{E}_x\left[f(B^+ \circ (t + \eta \gamma_t)^{-1})\right].
		\end{align*}
		Moreover, $\bar{e}_i$, $i \in \mathbb{N}$ are i.i.d. random variables for which (\cite[Section 4.3]{d2022fractionalsticky})
		\begin{align*}
			\mathbf{P}_0(\bar{e}_i >t | \bar{X}_{\bar{e}_i}>0) = \mathbf{E}_0\left[\exp\left(-(1/\eta) L^\Psi_t \right)\right], \quad \eta>0
		\end{align*}
		and
		\begin{align}
			\mathbf{E}_0[H^\Psi_{e_i}|B^+] = \mathbf{E}_0[e_1|B^+] \mathbf{E}_0[H^\Psi_1] = \mathbf{E}_0[e_1|B^+]\, \lim_{\lambda \downarrow 0} \frac{\Psi(\lambda)}{\lambda}
			\label{meanHoldingTime}
		\end{align}   
		gives the mean holding time on the boundary point $x=0$ for the process $\bar{X}$. That is, 
		\begin{align*}
			\upsilon(t,x) = \mathbf{E}_x\left[f(\bar{X}_t)\right] = \mathbf{E}_x\left[f(B^+ \circ (t + H^\Psi \circ \eta \gamma_t)^{-1}) \right],
		\end{align*}
		solves
		\begin{align*}
			\begin{cases}
				\displaystyle \dot{\upsilon}(t,x)= \upsilon^{\prime \prime}(t,x) & t>0,\, x \in (0,\infty), \\
				\displaystyle \eta\, \mathfrak{D}^\Psi_t \upsilon(t,x) = \upsilon^\prime(t,x) & t>0,\, x=0,\\
				\displaystyle \upsilon(0,x)=f(x), \; f \in C(0, \infty) \cap L^\infty(0, \infty) & x \in [0, \infty),
			\end{cases}
		\end{align*}
		for $\eta \geq 0$.\\
		
		Let us consider $\mathfrak{r}_t:=\sup\{s\geq 0\,:\, \gamma^+_s \leq t\}$ such that $\gamma^+ \circ \mathfrak{r}_t = t$. For a given $t$, $\mathfrak{r}_t$ is a stopping time for $B^+$. As in \cite{HSU}, we introduce the boundary process $\mathfrak{X}_t := B^+ \circ \mathfrak{r}_t$. That is the boundary trace process of $B^+$. The set of times \text{$J=\{t \geq 0\,:\, \mathfrak{X}_{t-} \neq \mathfrak{X}_t\}$} gives the times in which $\mathfrak{X}_t$ has a jump. The local time $\gamma^+_t$ is constant only for the excursions of $B^+$ on $(0,\infty)$, a jump for $\mathfrak{X}_t$ can be realized in $(\mathfrak{r}_{t-}, \mathfrak{r}_t)$ in case of non-empty set. Here $\mathfrak{X}_t$ is a pure jump process. We write $J_t = J \cap (0,t]$ and recall that a.s. $J$ is a dense countable subset of $(0,\infty)$. This can be associated with the countable jumps of a Cauchy process. The processes $X$ and $\bar{X}$ move along the path of $B^+$ for which the zero set $\{0\leq t < \infty :\, B^+_t=0\}$ has no isolated points. Assume $T^+ \in J_t$, then $B^+_{T^+}=0$. Moreover,
		\begin{align*}
			\sum_{ j \in J_t} e_j \quad \textrm{and} \quad \sum_{j \in J_t} \bar{e}_j
		\end{align*}
		respectively give the time the processes $X$ and $\bar{X}$ spend on the boundary point $\{0\}$ up to time $t$. \\
		
		We now define
		\begin{align*}
			B^\star_t := B^\bullet \circ S^{-1}_t, \quad t\geq 0
		\end{align*}
		and recall that also $B^\bullet$ can be identified as a process moving along the path of $B^+$. The new clock $S_t^{-1}$ is associated to a non-local dynamic boundary condition and it  introduces a sequence of holding times for $B^\star$. Since $B^\bullet$ jumps according with the jumps of $H^\Phi$, then once it reaches the boundary point $\{0\}$, it starts afresh after a jump. However, $B^\bullet$ maintains its Markovian nature (see \cite[Section 13]{ito1963brownian}). This does not hold for $B^\star$.\\

		Let us denote by $\{e^\star_i,\, i \in \mathbb{N}\}$ the sequence of holding times for $B^\star$.

		\begin{theorem}
			The holding times $e^\star_i$ are i.i.d. random variables for which
			\begin{align}
				\label{defHTstar}
				\mathbf{P}_{\mathfrak{s}}(e^\star_i > t | B^\star_{e^\star_i} \neq \mathfrak{s}) = \mathbf{E}_0\left[\exp\left(-(1/\eta) L^\Psi_t \right)\right], \quad \eta>0, \; \mathfrak{s} \neq 0.
			\end{align}
			In particular, $e^\star_1 \stackrel{d}{=} \bar{e}_1$.
			\label{thm:HTstar}
		\end{theorem}
		\begin{proof}
			The process $L^\Phi_{\gamma_t}:= L^\Phi \circ \gamma_t$ can be regarded as the local time of $B^\bullet$ at zero \cite[section 14]{ito1963brownian} and
			\begin{align*}
				B^\star_t \stackrel{law}{=} B^\bullet \circ (t + H^\Psi \circ \eta \gamma^\bullet_t )^{-1}  \quad \textrm{where} \quad \gamma^\bullet_t = \gamma_t(B^\bullet).
			\end{align*}
			The right-hand side of \eqref{defHTstar} says that $e^\star_i$ is an holding time at the point $\mathfrak{s} \in (0, \infty)$ and the left-hand side of \eqref{defHTstar} can be obtained by following the same arguments as in \cite[Lemma 6]{d2022fractionalsticky}. We say that $\mathfrak{s}$ is a sticky point for $B^\star$. For $\Psi=Id$, we get
			\begin{align*}
				\mathbf{P}_\mathfrak{s}(e^\star_i > t | B^\star_{e^\star_i} \neq \mathfrak{s}) = \mathbf{P}_\mathfrak{s}(e_i > t | B^\bullet_{e_i} \neq \mathfrak{s}) = \exp\left(-(1/\eta) t \right), \quad \eta>0, \; \mathfrak{s}\neq 0
			\end{align*}
			and we say that $\mathfrak{s} \in (0, \infty)$ is a sticky point for $B^\bullet$. Thus, the processes $B^\star$ and $B^\bullet$ near the origin jump to a point $\mathfrak{s}\neq 0$ and stop there according to the distribution of the associated holding time: $e_i$ are (i.i.d.) exponential holding times with parameter $(1/\eta)$ and $e^\star_i$ equal in law $H^\Psi_{e_i}$.
		\end{proof}
		
		We stress the fact that  $\mathfrak{s} \in B^\bullet \circ J$ and $B^\bullet \circ J$ can be regarded as the set of regenerative points for $B^\star$, that is we have
		\begin{align}
			\label{MarkovPropReg}
			\mathbf{E}_x\left[ f(B^\star_{t + J}) \big| B^\star_{J} \right] = \mathbf{E}_{B^\star_{J}} \left[ f(B^\star_t) \right], \quad t>0.
		\end{align}
		Let us discuss on this. The jumps of $B^\bullet$ are given by $\mathsf{L}_{\gamma_t}:= \mathsf{L} \circ \gamma_t$ as previously described. It may have jumps only as $B^+$ approaches $x=0$ and then $\gamma_t$ increases. Moreover, we know that $\mathbf{P}_0(\mathsf{L}_t > 0)=1$ for every $t>0$ and therefore the process $B^\star$ leaves the origin only through jumps. Moreover, the (length of the) jumps are independent from the base process $B^+$. The jumps determine the sticky points. Let us introduce the regenerative sets
		\begin{align*}
			\mathbf{T}_y = \{ (t, \omega) \subset \mathbb{R} \times \Omega \, : \, y=\mathsf{L}_{\gamma_{t}}(\omega) \}, \quad y \in (0, \infty)
		\end{align*}   
		and
		\begin{align*}
			\mathbf{S} = \{ (y, \omega) \subset \mathbb{R} \times \Omega \, : \, y=\mathsf{L}_{\gamma_{z}}(\omega) \textrm{ for some } z \geq 0 \}.
		\end{align*}   
		The random set $\{\mathbf{T}_y,\, y \in \mathbf{S}\}$ is not empty and includes the set $J$. In particular, the set
		\begin{align*}
			\{\mathsf{L}_\gamma \circ \tau\}, \quad \tau \in J
		\end{align*}
		gives the set of \textit{sticky points} for $B^\star$. Furthermore, $\{\tau \in J\}$ is a Markov set of points for $B^\bullet$ by means of which we may recover the associated Markov set of points for $B^\star$ by shifting w.r.t. holding times $\{e^\star_i -, \, : i \in \mathbb{N}\}$ and for which \eqref{MarkovPropReg} holds true. 
		The process $B^\star$ starts afresh at its Markov points. We observe that the set of sticky points has zero Lebesgue measure whereas, the set of holding times is a set of positive Lebesgue measure (which obviously does not hold for the set $\{t\,:\, B^+_t=0\}$).\\

		In conclusion, the solution to the NLBVP in Theorem \ref{fracdyntm} has a probabilistic representation given by $B^\star$ which can be constructed via time-change $S$ of a decomposition in terms of $B^+$ and $\mathsf{L}_\gamma$. The process moves along the path of  the reflected Brownian motion. The time-change gives the amplitude of the plateaus, it is governed by the non-local operator in time. The additive part gives the amplitude of the jumps, it is governed by the non-local operator in space. Reflected Brownian motion, jumps and plateaus are independent.

		\subsection{Stochastic resetting}
		\label{sec:resetting}
		In this section we underline that the time reverse process of $B^\bullet$ produces paths which can be associated with Brownian motions under stochastic resetting. 
		\begin{figure}
			\centering
			\begin{minipage}{0.5\linewidth}
				\centering
				\includegraphics[width=.9\linewidth]{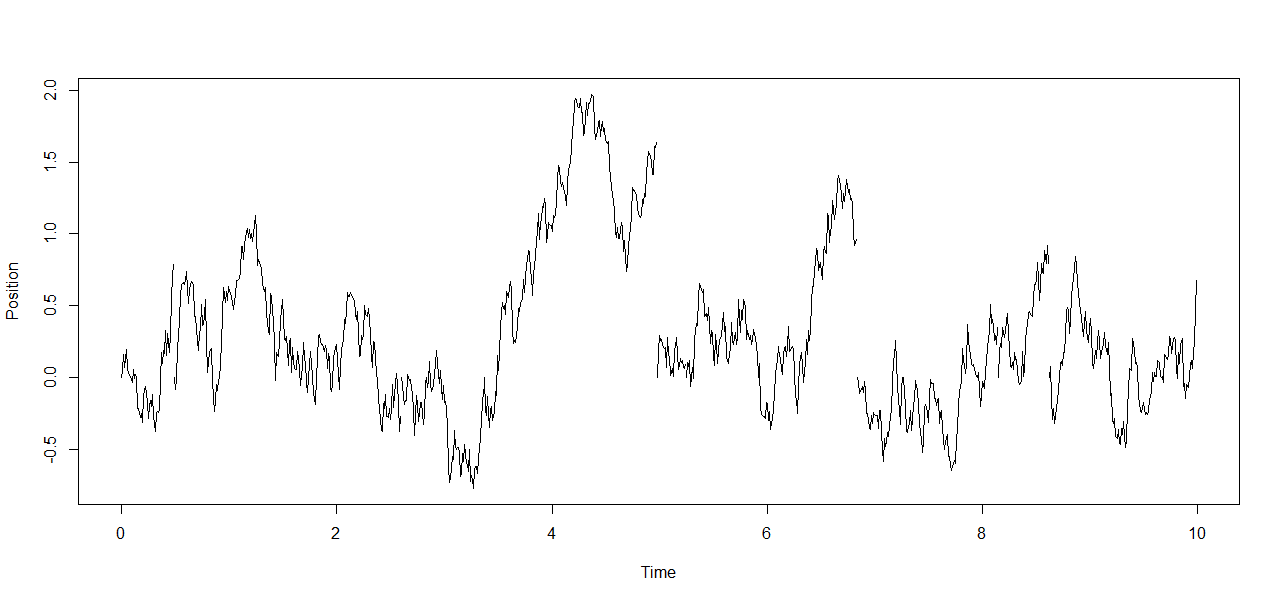}
			\end{minipage}\hfill
			\begin{minipage}{0.5\linewidth}
				\centering
				\includegraphics[width=.9\linewidth]{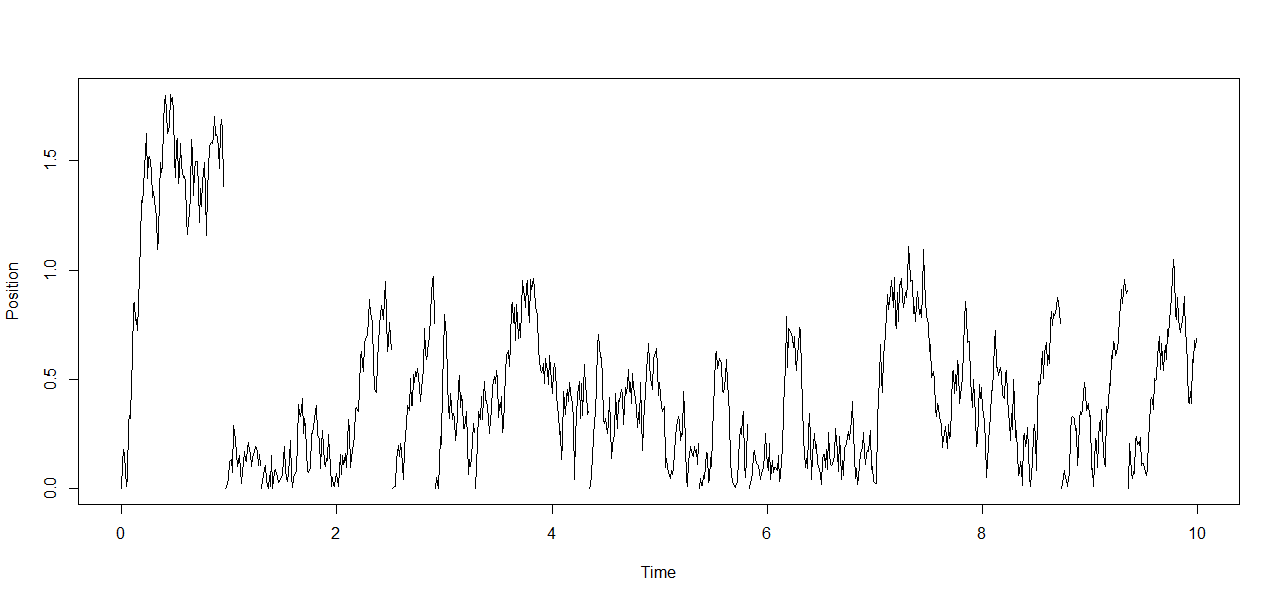}
			\end{minipage}
			\caption{Two examples of independent Brownian motions with Poissonian resetting. The left plot shows the process on $\mathbb{R}$, while the right plot illustrates the process constrained to the positive half-line.}
			
			\label{resetting}
		\end{figure}
		
\begin{figure}
			\centering
			\begin{minipage}{0.5\linewidth}
				\centering
				\includegraphics[width=.9\linewidth]{BMtimechangenewnew.png}
			\end{minipage}\hfill
			\begin{minipage}{0.5\linewidth}
				\centering
				\includegraphics[width=.9\linewidth]{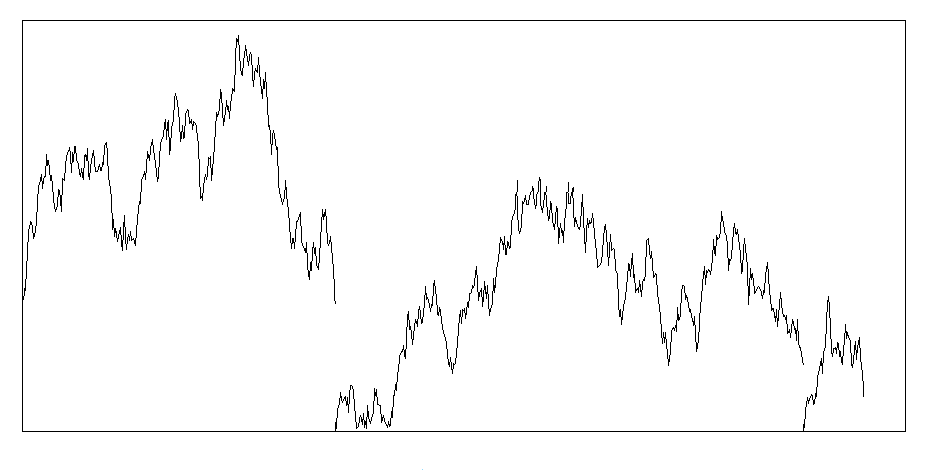}
			\end{minipage}
			\caption{A possible path for $B^\bullet_t$  $0 \leq t \leq T$ (left) and $B^\bullet_{T-t}$ , $0\leq t \leq T$ reverse over time (right). From the path in Figure \ref{Bpallino} we obtain the paths which has been drown here.}
			\label{BpallinoSR}
		\end{figure}

		Let us consider the Figure \ref{Bpallino} given by the path of $B^\bullet_t$ with $0\leq t \leq T$ for a given finite time $T>0$. We introduce the process $\mathfrak{B}=\{\mathfrak{B}_t\}_{0 \leq t \leq T}$ defined as
		\begin{align*}
			\mathfrak{B}_t := B^\bullet_{T-t}, \quad 0 \leq t \leq T
		\end{align*}
		where $\mathfrak{B}_0 = B^\bullet_T$ and $\mathfrak{B}_T=B^\bullet_0 = x>0$. This produces the paths in Figure \ref{BpallinoSR}. Let \text{$J=\{J_t\}_{0\leq t \leq T}$} give the random jumps of $B^\bullet$ from zero. For example \text{$J=\{t\,:\, B^\bullet_{t-} \neq B^\bullet_t\}$}. Recall that $\mathfrak{B}$ behaves like $B^\bullet$ on $(0, \infty)$. Moreover, in both cases the excursions on $(0, \infty)$ are Markovian, that is $B^\bullet$ is Markov (\cite[Section 13]{ito1963brownian}) as well as $\mathfrak{B}$. We focus on the pictures of Figure \ref{BpallinoSR}. The associated process $B^\bullet$ (started at $B^\bullet_0=x>0$) behaves like a Brownian motion on $(0, \infty)$ and
		\begin{align*}
			\inf\{t \,:\, B^\bullet_t=0\} =:\tau_0
		\end{align*} 
		has the distribution $\mathbf{P}_x(\tau_0 >t) = Q^D_t \mathbf{1}(x)$ where $Q^D_t$ has been introduced in \eqref{upsilonSol}. The zero set of a Brownian motion is a Cantor set. However, because of the jumps away from zero, $\tau_0$ and in general, the sequence $\{\tau_i\}_{ i \in \mathbb{N}}$ of (zero) hitting times are well identified. For the process $B^\bullet$ we can say that $\tau_i \sim \tau_0$ $\forall\, i$ (the times are i.i.d.).  The process $\mathfrak{B}$ (started at $\mathfrak{B}_0=B^\bullet_T$) behaves like the Brownian motion on $(0, \infty)$ until a resetting occurs. As suggested by the right picture in Figure \ref{BpallinoSR}, the resetting of $\mathfrak{B}$ can be associated with the set $J$ of the jumping times of $B^\bullet$. Since the resetting level is given by the jumps of $B^\bullet$ with each visit of the boundary point $\{0\}$, we have that
		\begin{align}
			\label{Hy}
			\mathbf{P}(B^\bullet \textrm{ has a jump away from zero at time }  > t \, |\,  B^\bullet_{0} = \mathfrak{s}) = \mathbf{P}_{\mathfrak{s}}(\tau_0 >t), \quad \mathfrak{s} \in B^\bullet \circ J
		\end{align}
		coincides with
		\begin{align}
			\label{Hx}
			\mathbf{P}(\textrm{ the resetting time $\mathfrak{T}$ of $\mathfrak{B}$ to zero }  > t\, \, |\, \mathfrak{B}_\mathfrak{T}=\mathfrak{s}) = \mathbf{P}_{\mathfrak{s}}(\tau_0 >t), \quad \mathfrak{s} \in B^\bullet \circ J
		\end{align}
		where $\mathfrak{T}$ is the resetting time random variable taking values in $T- J$. In particular, we recall that  
		\begin{align*}
			\mathbf{P}_x(\tau_0>t) = \mathbf{P}_0(H_x >t) = \mathbf{P}_0(x < L_t), \quad x \in (0, \infty)
		\end{align*} 
		where $H_x$ is a stable subordinator of order $\alpha=1/2$ and $H_x \stackrel{d}{=} \inf\{s\,:\, B_s=x\}$. We underline that $B^\bullet \circ J$ is a random set obtained from the path of $B^\bullet$ at times given in $J$, that is the jumping times of $B^\bullet$. Thus, \eqref{Hy} and \eqref{Hx} can be considered for a r.v. $\mathfrak{s}>0$. Notice that we only have equivalence in law for the hitting times above, the underlying paths of the associated Brownian motions are reasonably different.
		Let us consider the case $\Pi^\Phi(0,\infty)<\infty$, so that the subordinator $H^\Phi$ is not of pure-jump type (in particular $d>0$ in \eqref{LevKinFormula}). 
		From the definition \eqref{process} it follows that, at zero, the process $B^\bullet$ may either reflect or jump: it jumps whenever the subordinator has a jump and it is reflected otherwise. 
		This behaviour is consistent with the picture in Figure \ref{resetting}: with Poissonian resetting the forward process can hit zero either by ordinary Brownian motion or as the result of a reset. 
		Consequently, in the time-reversed dynamics one also observes reflection on those sample paths where the forward process reached zero without a reset. 
		\\
		In fact, it was shown in \cite{ColDovPag} that the time-reversed process is characterized by a NLBVP that incorporates both reflection and jumps. 
		Moreover, the backward dynamics feature an additional negative drift term: this drift compensates for the concentrated mass of the reset process near zero.
		
		\section{Motions on graphs}
		\label{sec:graphs}
		
		Brownian motions on metric graphs describe the movement of particles across interconnected nodes. We consider non-local conditions on nodes, which is the novel here,  by first studying a fundamental case with a single node, that is a star graph.
		The star graph is a special metric graph where there is a finite collection of sets, isomorphic to the positive half-line, such that all the origin points coincide at the unique point which is the vertex. Since each edge is isomorphic to $(0,\infty)$, it is natural to extend Feller's theory on boundary conditions \cite{feller1952parabolic} and the related probabilistic solutions \cite{ito1963brownian}. The pathwise construction of the Brownian motion on the star graph and the analysis of conditions at the vertex are developed in \cite{russi1,russi2}. The topic has been thoroughly covered also in the thesis \cite{werner2016brownian}. The construction for a general metric graph is treated in \cite{russi3}. Recently, the sticky Brownian motion, even with non-local dynamic conditions, has been studied in \cite{BonDov}.
		
		\subsection{Star graph and Brownian motions}
		\label{setting}
		Let us begin by introducing the concept of a star graph within the theory of metric graphs. Let $\mathcal{E}$ be a countable finite set and ${E}$ be a family of copies of the positive half-line
		\begin{align*}
			{E}:= \bigsqcup_{e \in \mathcal{E}} [0,\infty).
		\end{align*}
		For simplicity, we denote each point in ${E}$ by a couple $(j,x)$, where $j \in \mathcal{E}$ and $x$ is the Euclidean distance from the origin. As in \cite{mugnolo2019actually}, we introduce the following equivalence relation
		\begin{align*}
			(j,x) \sim (i, y) \iff \begin{cases}
				j=i \ \text{and } x=y\\
				x=y=0, \ \text{for any } j,i.
			\end{cases}
		\end{align*}
		Then, the star graph is the quotient space $\mathcal{G}:= {E}/ \sim$. The unique vertex $v=(\cdot,0)$ belongs to all the edges. We also have that $\mathcal{G}$ is a metric space with the distance
		\begin{align*}
			d((j,x),(i,y)):=\vert x-y \vert \mathbf{1}_{j=i} + (x+y) \mathbf{1}_{j \ne i}.
		\end{align*}
		Another possible representation for the star graph $\mathcal{G}$, with star vertex $v$, is the following
		\begin{align*}
			\mathcal{G}=\{v\} \cup \bigcup_{e \in \mathcal{E}} (\{e\} \times (0, \infty) ).
		\end{align*}
		\begin{figure}
			\centering
			\includegraphics[width=8cm]{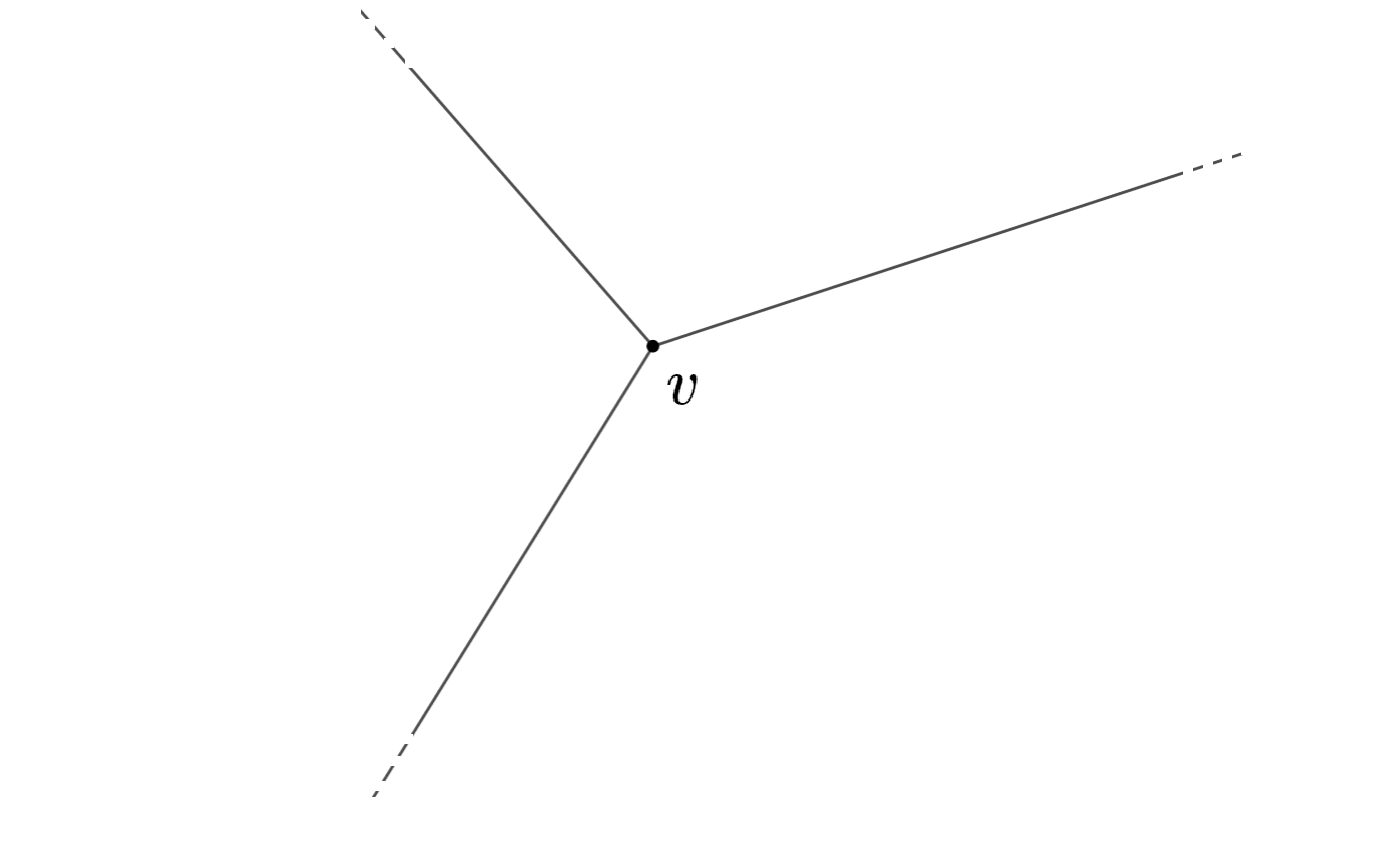}\\
			\caption{A star graph with $3$ edges and vertex $v$.}
		\end{figure}
		Now, we define the Brownian motion on the star graph $\mathcal{G}$. With $\epsilon_e$, we indicate
		\begin{align*}
			\epsilon_e(\mathsf{x})=\epsilon_e(j,x)=\begin{cases}
				1 \quad &\text{if } j=e\\
				0 \quad &\text{otherwise}
			\end{cases}.
		\end{align*}
		
		\begin{theorem}
			\label{thm:equiv}
			({\it Equivalence}) Let $\mathcal{T}^U:=\inf \{t>0 : \mathcal{B}_t=v\}$ be the first hitting time for the Brownian motion $\mathcal{B}$ on $\mathcal{G}$ of the vertex $v$, that is the point $(\varepsilon, 0) \equiv 0 \in \mathcal{G}$ and write $\mathcal{T}_0 = \mathcal{T}_v$. The process $\mathcal{B} \circ (t \wedge \mathcal{T}_0)$ started at $(\varepsilon, x)$ is equivalent in law to the process $B^+ \circ (t \wedge \tau_0)$ started at $x$ for any $x \neq 0$ and $\varepsilon \in \mathcal{E}$ with $\tau_0 = \inf\{t>0\,:\, B^+_t = 0\}$.
		\end{theorem}
		\begin{proof}
			The proof follows from the same arguments as in Theorem 4 in \cite{BonDov}. Notice that both processes behave like a Brownian motion away from $0$.
		\end{proof}
		
		\begin{definition}(Brownian motion on a star graph)
			\label{bm:graph}
			Let $\mathcal{B}=(U, B^+)$ be defined on $\mathcal{G}$, such that 
			\begin{align*}
				\mu_U := \sum_{e \in \mathcal{E}} \frac{1}{\vert \mathcal{E} \vert} \epsilon_e.
			\end{align*}
			The process $\mathcal{B}$ satisfies:
			\begin{itemize}
				\item $B^+$ is a reflected Brownian motion;
				\item If $\mathcal{B}_0=v$, then for $t>0$, the distribution of $U_t$ is given by $\mu_U$;
				\item For $\mathcal{B}_0=(l,x)$, for $x>0$, then  $U_t=l$ if $t \leq \mathcal{T}^U$ and if $t>\mathcal{T}^U$ the distribution of $U$ is equal to $\mu_U$ and independent of $B^+$.
			\end{itemize}
		\end{definition}
		\begin{remark}
			\label{walsh}
			If in Definition \ref{bm:graph}, we consider $\mathcal{X}=(\Theta,B^+)$ and $\Theta$ is not a uniform random variable on the external edges, for example if $\Theta \sim \mu$ such that $\mu := \sum_{e \in \mathcal{E}} \rho_e \epsilon_e$ and $\sum_{e \in \mathcal{E}} \rho_e = 1$, then we are dealing with a Walsh Brownian motion on the star graph \cite{walsh1978diffusion}.
		\end{remark}
		\subsection{Functions on star graphs}
		Contrary to Brownian motion on the real line, where we can visualize a particle free to move, on the star graph, the presence of the vertex introduces necessarily boundary conditions. First of all, we must grasp the concepts of functions, continuity, and differentiability on the graphs.
		
		Let us introduce the spaces of functions that we need. Since $(\mathcal{G},d)$ is a metric space, a function $f: \mathcal{G} \to \mathbb{C}$ is continuous if the standard definition of continuous functions on metric spaces holds. In addition, we observe that $f$ can be considered as
		\begin{align*}
			f=\bigoplus_{e \in \mathcal{E}} f_e,
		\end{align*}
		where $f_e(\cdot)= f(e,\cdot) : (0,\infty) \to \mathbb{C}$ is the projection of $f$ on the edge $e$, for \(v=(e,0)\), we define \(f(v)=f_e(0)\). Then for all $k \in \mathbb{N}$ we define recursively
		\begin{align*}
			\begin{split}
				&C(\mathcal{G}):=\left\{f=\bigoplus_{e \in \mathcal{E}} f_e \in \bigoplus_{e \in \mathcal{E}} C(0,\infty) \text{ and } f_e(0)=f_l(0), \, \forall e, l \in \mathcal{E}\right\}\\
				&C^k(\mathcal{G}):=\left\{f=\bigoplus_{e \in \mathcal{E}} f_e \in \bigoplus_{e \in \mathcal{E}} C^k(0,\infty) \text{ and } f_e^{(k)}(0)=f_l^{(k)}(0), \, \forall e, l \in \mathcal{E}:\right.
				\\
				&\left. f^{(h)} := \bigoplus_{e \in \mathcal{E}} f^{(h)}_e \in C(\mathcal{G}), 1 \leq h \leq k \right\},
			\end{split}
		\end{align*}
		with $f^{(h)}_e$ we denote the $h-$derivative of $f_e$. In this way the functions are continuous up to the vertex. 

		Similarly to what we observed for functions, measures can also be regarded as direct sums of measures on edges. Consequently, we can, for instance, define Lebesgue spaces. By following \cite[Definition 3.15]{mugnolo-libro}, we denote by $L^p(\mathcal{G})$ the space of measurable functions $f: \mathcal{G} \to \mathbb{C}$ such that
		\begin{align*}
			&\vert \vert f \vert \vert_{L^p(\mathcal{G})} :=\left(\sum_{e\in \mathcal{E}} \int_0^\infty \vert f(e,x)\vert^p dx\right)^\frac{1}{p}<\infty \quad p\in [1,\infty)\\
			&\vert \vert f \vert \vert_{L^\infty(\mathcal{G})}:=\inf \{c \in \mathbb{R}: \vert f(e,x)\vert \leq c \text{ for a.e. } x\in (0,\infty), \text{ and for all }e \in \mathcal{E}  \}.
		\end{align*}
		Now, the Sobolev space is, for $ p\in [1,\infty]$,
		\begin{align*}
			W^{k,p}(\mathcal{G})&:=\left\{f \in L^p(\mathcal{G}), f^{(h)} \in C(\mathcal{G}) \text{ for } 1 \leq h \leq k-1,\right.\\
			&\left. \text{ and } f^{(j)} \in L^p(\mathcal{G}) \text{ for all } 0 \leq j \leq k \right\}.
		\end{align*}
		Let us turn our attention to the concept of derivatives. In our exploration of the heat equation on graphs, we direct our focus toward the Laplacian.
		Let $C_0(\mathcal{G})$ be the space of the continuous functions $f$ such that, for every $e \in \mathcal{E}$, $f_e$ vanishes at infinity. The Laplacian is the operator $\Delta: C^2_0(\mathcal{G}) \to C_0(\mathcal{G})$ written as
		\begin{align}
			\label{laplacian}
			\Delta u (\mathsf{x})= \frac{d^2}{dx^2} u(l,x)
		\end{align}
		where $\mathsf{x}=(l,x) \in \mathcal{G}$ and $\mathsf{x}\ne v$. For the star vertex $v=(e,0)$, with  $e \in \mathcal{E}$
		\begin{align}
			\label{laplacian-vertex}
			\Delta f(v):=\lim_{\xi \to 0} f^{\prime \prime}_e(\xi),
		\end{align}
		for which the limit exists and it is a function in $C_0(\mathcal{G})$. Representing the vertex as $v=(e,0)$, when computing the derivative, it is crucial to keep track of the direction in which we are moving. For this reason, we take the limit along the projection onto the edge $e$. In the case of the Laplacian, we stipulate that it is a continuous operator across the entire graph, ensuring that the limit \eqref{laplacian-vertex} holds the same value for every edge $e \in \mathcal{E}$. Similarly, we can define other derivatives, as we have done for the second one, and require them to be continuous outside the vertex; in such instances, the limits, like \eqref{laplacian-vertex}, depend on the edge $e$. For example, for $u \in C^1(\mathcal{G} \setminus \{v\})$
		\begin{align*}
			u^\prime (\mathsf{x}):=\frac{d}{dx} u(l,x)
		\end{align*}
		with $\mathsf{x}=(l,x)$. 
		But, in this case, we are not considering that at the vertex, the derivatives are all equal.
		\subsection{Stochastic processes on the edges}
\begin{figure}[!b]
			\centering
			\includegraphics[width=8cm]{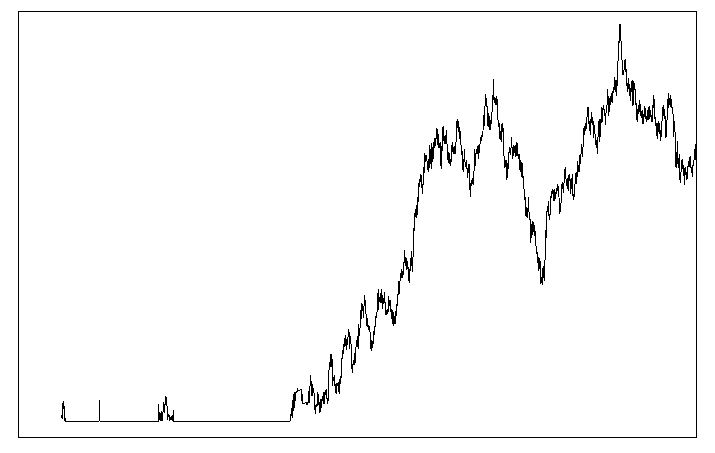} 
			\caption{A possible path for $X^e$.}
			\label{paths-X^e}
		\end{figure}

		In Definition \ref{bm:graph}, our focus was on Brownian motion, examining it across the entire graph. Presently, we aim to comprehend its behavior on each individual edge.\\
		Formally, we define the process $X_t^e$, on the edge $e \in \mathcal{E}$ as
		\begin{align*}
			X^e_t= f(\mathcal{B}_t),
		\end{align*}
		where $f : \mathcal{G} \to \mathbb{R}$ is
		\begin{align*}
			f(\mathcal{B}_t):=\begin{cases}
				B^+_t \quad &\text{if } U_t=e\\
				0 \quad &\text{if } U_t\ne e.
			\end{cases}
		\end{align*}
		Based on the definition, it is evident that $X^e$ behaves as a reflecting Brownian motion up to the point where $U=e$. However, when the Brownian motion $\mathcal{B}$ shifts between edges on the graph, the $X^e$ process remains stuck at zero until $\mathcal{B}$ occurs again on the same edge $e$. As illustrated in Figure \ref{paths-X^e}, we observe that the plateaus (intervals of consistency) of $X^e$ are a result of excursions of $\mathcal{B}$ to different edges. This observation leads us to anticipate a non-Markovian dynamic.
		We observe that this case is very close to the case \( \eta > 0 \) and \( \Phi = \text{Id} \) already studied (see Figure \ref{onlysticky}). In both cases, we deal with processes that remains stuck at zero: for \( X^e \), the slowdown is caused by the excursions of \( \mathcal{B} \) on the other edges, while for \( B^+ \circ \bar{V}^{-1} \) the slowdown is due the non-local dynamic boundary conditions.
		\begin{theorem}
			Then the process $X^e$ is not Markov with respect to its natural filtration.
		\end{theorem}
		\begin{proof}
			Fix a time $t\ge 0$ and define the random time
			\begin{align*}
				T^e_t:=\inf\{s>t:\;U_s=e\},
			\end{align*}
			and the event
			\begin{align*}
				A:=\{T^e_t-t=0\}=\{\inf\{s>t:\;U_s=e\}=t\}.
			\end{align*}
			Consider the two events
			\begin{align*}
				H_1:=\{U_t=e,\;B^+_t=0\},\qquad H_2:=\{U_t\neq e,\;B^+_t>0\}.
			\end{align*}
			Both $H_1$ and $H_2$ imply the same observed value $X^e_t=0$, because $X^e_t$ equals $0$ when the process is on another edge and also when it is at the vertex with distance zero.\\
			However, the probability of $A$ is different under $H_1$ and $H_2$:
			if $H_1$ holds, then the process is already on edge $e$ at time $t$, so $T^e_t=t$ almost surely and
			\begin{align*}
				\mathbf{P}(A\mid H_1)=1.
			\end{align*}
			If $H_2$ holds, then the process is on a different edge at positive distance from the vertex, so it must first return to the vertex and then choose an edge; hence $T^e_t>t$ almost surely and
			\begin{align*}
				\mathbf{P}(A\mid H_2)=0.
			\end{align*}
			If $X^e$ was Markov with respect to its natural filtration, the conditional law of any event (in particular $A$) given the past up to time $t$ would depend only on the current value $X^e_t$. Thus we would have
			\begin{align*}
				1= \mathbf{P}(A\mid H_1)=\mathbf{P}(A\mid X^e_t=0)=\mathbf{P}(A\mid H_2)=0,
			\end{align*}
			a contradiction. Therefore $X^e$ is not Markov.
		\end{proof}

		\subsection{Fundamental solution of the heat equation}
		In this section, we analyze the heat equation on the star graph and observe how only with the choice of uniform weights on the edges we achieve the expected inclusion in  operator's domain.
		
		Let $P_t^{\mathcal{X}} f(\mathsf{x})= \mathbf{E}_{\mathsf{x}}[f(\mathcal{X}_t)]$ be the Feller semigroup of the Walsh Brownian motion $\mathcal{X}$, see \cite[Theorem 2.1]{yor-walsh}, defined in Remark \ref{walsh}, for $f \in C(\mathcal{G})$ and $\mathsf{x}=(l, x)$. From \cite[Formula (2.1) and (2.2)]{yor-walsh}, we have
		\begin{align*}
			\displaystyle
			P_t^\mathcal{X} f(\mathsf{x})= \sum_{e\in \mathcal{E}} \rho_e \left(P_t^+ f_e(x) - P_t^D(f_l  - f_e)(x)\right),
		\end{align*}
		where $P^+_t$ is the semigroup of a reflected Brownian motion on the half-line and $P_t^D$ is the semigroup of the Brownian motion on the half-line killed in zero. Then, the infinitesimal generator of $\mathcal{X}$ is $A=\Delta$, with $\Delta$ the laplacian on the star graph defined in \eqref{laplacian}, and the domain is
		\begin{align}
			\label{domain-walsh}
			D(A)=\left\{\varphi, \Delta \varphi \in C_0(\mathcal{G}) \ : \ \sum_{e \in \mathcal{E}} \rho_e \frac{\partial}{\partial x} \varphi(e,x)\big\vert_{x=0}=0\right\},
		\end{align}
		so that the Kirchoff's condition on the vertex appears.
		\begin{theorem}
			The density of the Walsh Brownian motion on $\mathcal{G}$ is not a function in $D(A)$.
		\end{theorem}
		\begin{proof}
			For the point $\mathsf{x}=(l,x)$, the density of the Walsh Brownian motion is given by
			\begin{align*}
				g(t,v, \mathsf{x})\, d\mathsf{x} :=\mathbf{P}(\mathcal{X}_t \in d\mathsf{x} \vert \mathcal{X}_0=v)= \rho_l \mathbf{P}(B^+_t \in dx \vert B^+_0=0)= 2 \rho_l g(t,x)\, dx,
			\end{align*}
			where $g$ is the Gaussian kernel. Therefore, we show that the heat equation holds
			\begin{align*}
				\frac{\partial}{\partial t} g(t,v, \mathsf{x})= 2 \rho_l \frac{\partial}{\partial t} g(t,x) = 2 \rho_l \frac{\partial^2}{\partial x^2} g(t,x)= \Delta g(t,v, \mathsf{x}),
			\end{align*}
			but if we want continuity at the vertex
			\begin{align}
				\label{lim_vertex}
				\lim_{\mathsf{x}\downarrow v} g(t, v, \mathsf{x})=\lim_{x \to 0}  2 \rho_l g(t,x) \text{ depends on the edge $l$} 
			\end{align}
			then it is not continuous and $g(t,v, \mathsf{x}) \notin C(\mathcal{G})$. So the density of the Walsh Brownian motion does not belong to $D(A)$.
		\end{proof}
		
		On the contrary, when we consider the Brownian motion $\mathcal{B}$, defined in Definition \ref{bm:graph}, the limit \eqref{lim_vertex} does not depend on the edge, since the distribution of edge selection is uniform. For the other derivatives, we provide
		\begin{align*}
			\lim_{\mathsf{x}\downarrow v} \frac{d}{dx} g(t,v, \mathsf{x})= \lim_{x \downarrow 0} \frac{\partial}{\partial x}g(t,v,(l,x))=\lim_{x \downarrow 0} \frac{2}{\vert \mathcal{E} \vert} (-g(t,x))=0
		\end{align*}
		and
		\begin{align*}
			\Delta g(t, v, \mathsf{x}) = \frac{2}{\vert \mathcal{E} \vert} (-g(t,x) + x^2 g(t,x)) \rightarrow - \frac{2}{\vert \mathcal{E} \vert} g(t,0) \quad {x \downarrow 0}
		\end{align*}
		which does not depend on the edge $l$. We conclude that $g(t,v,\mathsf{x}) \in D(A)$ if and only if we deal with a Brownian motion, such that the distribution of the edge selection is uniform.
		
		\subsection{Non-local operators on the star graph}
		Here we introduce the non-local operators on the star graph $\mathcal{G}$ as a natural extension of the results in Section \ref{sec:Marchaud} and Section \ref{sec:aux}.
		
		For $u \in W^{1,\infty}(\mathcal{G})$,  we define 
		the left  Marchaud-type derivative  on $\mathsf{x}=(e,x) \in \mathcal{G} \setminus \{v\}$
		\begin{align}
			\label{Marchaudleft-graphs}
			\mathbf{D}_{\mathsf{x}-}^\Phi u(g) = \int_0^\infty (u(e,x) - u(e,x+y))\Pi^\Phi(dy).
		\end{align}
		In this way, the function $u$ is $W^{1,\infty}$ on each edge and we use the analogue of \eqref{stimaMarchaud} to see that  \eqref{Marchaudleft-graphs} is well defined.
		
		As for the Marchaud (type) operator, we exploit the same idea for Riemann-Liouville (type) operators, which are
		\begin{align}
			\label{R-LxinfG}
			\mathcal{D}_{\mathsf{x}-}^\Phi u(x):=-\frac{d}{dx} \int_x^\infty u(e,y) \overline{\Pi}^\Phi(y-x) dy
		\end{align}
		and
		\begin{align}
			\label{R-L0xG}
			\mathcal{D}_{\mathsf{x}+}^\Phi u(g):=\frac{d}{dx} \int_0^x u(e,y) \overline{\Pi}^\Phi(x-y) dy
		\end{align}
		respectively defined for function $u$ such that 
		\begin{align*}
			u(e,\cdot) \overline{\Pi}^\Phi(\cdot-x) \in L^1(x,\infty), \quad \text{and} \quad u(e,\cdot) \overline{\Pi}^\Phi(x-\cdot) \in L^1(0,x) \quad \forall x \in (0,\infty), e \in \mathcal{E}.
		\end{align*}
		Regarding the star vertex $v$, careful consideration of the definitions is essential. We represent $v$ as $(e,0)$ across all edges, necessitating specification of the direction in which we are progressing. For example, for \eqref{Marchaudleft-graphs}, we have for $u \in W^{1,\infty}(\mathcal{G})$
		\begin{align}
			\label{marchaud-vertex}
			\mathbf{D}^\Phi_{v-} u(v) = \lim_{\xi \to 0} \mathbf{D}^\Phi_{\xi-} u(e,\xi)= \lim_{\xi \to 0} \mathbf{D}^\Phi_{\xi-} u_e(\xi),
		\end{align}
		with  $\mathbf{D}^\Phi_{\xi-} u(e,\xi)$ is the edgewise Marchaud-type derivative \eqref{Marchaudleft} on $\xi$.
		\subsection{Non-local conditions on star graphs}
		As previously anticipated, the fact that each edge is isomorphic to the positive half-line leads us to think that Feller's theorem is easily extendable to the star graph $\mathcal{G}$. This is true, see for example \cite[Theorem 20.26]{werner2016brownian}. In general, the integral boundary condition is typically interpreted as the integral over a measure on the graph $\mathcal{G}$. However, we prefer to focus on L\'evy measures and provide an interpretation involving non-local operators.
		\begin{definition}
			\label{def:Xpallino}
			Let $\mathcal{X}^\bullet=(\Theta, B^\bullet)$ be defined on $\mathcal{G}$, with weights $\rho_e \geq 0$ for all $e \in \mathcal{E}$ such that 
			\begin{align*}
				\sum_{e \in \mathcal{E}} \rho_e = 1 \text{ and } \mu := \sum_{e \in \mathcal{E}} \rho_e \epsilon_e.
			\end{align*}
			The process $\mathcal{X}^\bullet$ satisfies:
			\begin{itemize}
				\item $B^\bullet$ is the process defined in \eqref{process};
				\item If $\mathcal{X}^\bullet_0=v$, then for $t>0$, the distribution of $\Theta_t$ is given by $\mu$;
				\item For $\mathcal{X}^\bullet_0=(l,x)$, for $x>0$, then  $\Theta_t=l$ if $t \leq \mathcal{T}$ and if $t>\mathcal{T}$ the distribution of $\Theta$ is equal to $\mu$ and independent of $B^\bullet$,
				where 
				\begin{align*}
					\mathcal{T}:=\inf \{t>0 : \mathcal{X}^\bullet_t=v\}.
				\end{align*}
			\end{itemize}
		\end{definition}
		Since $B^\bullet$ is a  Markov process, the same property is inherited by $\mathcal{X}^\bullet$. $\mathcal{X}^\bullet$ starts its paths in a point $\mathsf{x}$, it behaves like a Walsh Brownian motion on a star graph and, when it reaches the vertex $v$, it jumps (as the last jumps of $H^\Phi$) in an edge (by choosing it trough the distribution of $\Theta$). On each edge the paths are right-continuous, since the same is true for $B^\bullet$.
		\begin{remark}
			This definition could be also obtained from the Walsh Brownian motion on the star graph, as in \cite[Section 21]{werner2016brownian}.
		\end{remark}
		\begin{remark}
			The second component of $\mathcal{X}^\bullet$ on $\mathcal{G} \setminus \{v\}$ and $B^\bullet$ on $(0,+\infty)$ have the same excursions. Indeed, from the definition of $\mathcal{X}^\bullet$, we have that outside from the node $v$ the process coincides with $B^\bullet$.
		\end{remark}
		We have seen that, even for simple Brownian motion on $\mathcal{G}$, it is necessary to impose Kirchhoff conditions at the vertex. Therefore, it makes sense to search solutions in the space
		\begin{align*}
			\mathcal{D}:=\left\{\varphi(\cdot,\mathsf{x}) \in C(0,\infty) \, \forall\mathsf{x}\in \mathcal{G}; \, \varphi(t, \cdot), \Delta\varphi(t, \cdot) \in C_0(\mathcal{G}), \, \varphi^\prime(t, \cdot) \in C(\mathcal{G}\setminus\{v\}) \, \forall t >0 \right\}.
		\end{align*}
		For the non-local Kirchoff condition, we need $D_H \subset \mathcal{D}$ defined as
		\begin{align*}
			D_H:=\left\{\varphi \in \mathcal{D}: \varphi(t,\cdot) \in W^{1,\infty}(\mathcal{G}), \, \mathbf{D}^\Phi_{x-} \varphi_e(t,x)\big\vert_{x=0} \, \text{exists} \, \forall e \in \mathcal{E}, \, \forall t>0 \right\},
		\end{align*}
		where the Marchaud (type) derivative at the vertex $v$ is rewritten trough the Marchaud (type) derivatives on the projections $\varphi_e$, as in \eqref{marchaud-vertex}.
		\begin{theorem}
			\label{thm:Xpallino}
			The probabilistic representation of the solution $u \in D_H$
			\begin{align}
				\label{problem}
				\begin{cases}
					\dot{u}(t,\mathsf{x})= \Delta u(t,\mathsf{x}) \quad &t>0, \mathsf{x} \in \mathcal{G}\setminus \{v\} \\
					\sum_{e \in \mathcal{E}} \rho_e\mathbf{D}^\Phi_{x-} u_e(t,x)\big\vert_{x=0}=0 \quad    &t>0\\
					u(0,\mathsf{x})=f(\mathsf{x}) \quad &\mathsf{x} \in \mathcal{G}
				\end{cases}
			\end{align}
			with $f \in C(\mathcal{G}) \cap L^\infty(\mathcal{G})$  and $\mathsf{x}=(l,x)$ is
			\begin{align*}
				u(t,\mathsf{x})=\mathbf{E}_\mathsf{x}[f(\mathcal{X}^\bullet_t)]
			\end{align*}
		\end{theorem}
		We postpone the proof in Section \ref{proof2}
		\begin{remark}
			In this remark, we delve into how our process can be viewed as a generalization of the Walsh Brownian motion on the star graph, by taking advantage of the properties of fractional calculus.
			Let $f$ be a function on $\mathcal{G}$ such that \text{$f \in W^{1,1}(\mathcal{G})$} and $\Phi(\lambda)=\lambda^\alpha$ with $\alpha \in (0,1)$. Then, the associated L\'evy measure, from \eqref{LevKinFormula}, is denoted by $\Pi^\alpha$ and $\mathbf{D}_{\mathsf{x}-}^\alpha f(\mathsf{x})$ is the Marchaud derivative, with $\mathsf{x}=(e,x) \in \mathcal{G}$. We see that
			\begin{align*}
				\mathbf{D}_{\mathsf{x}-}^\alpha f(\mathsf{x})&:=\int_{0}^{\infty} (f(e,x) - f(e,x+y))\Pi^\alpha(dy)\\
				&=- \int_{0}^{\infty} \frac{d}{dz} f(e,x+z) dz \, \Pi^\alpha(dy)\\
				&=-\int_{0}^{\infty} \frac{d}{dz} f(e,x+z) \Pi^\alpha((z,\infty)) dz\\
				&=\text{[by using \eqref{tailSymb}]}\\
				&=-\int_{0}^{\infty} \frac{d}{dz} f(e,x+z) \overline{\Pi}^\alpha (z) dz\\
				&=-\int_{x}^{\infty} \frac{d}{dz} f(e,z) \overline{\Pi}^\alpha(z-x) dz.
			\end{align*}
			We observe that
			\begin{align*}
				\int_{0}^{\infty} e^{-\lambda z} \overline{\Pi}^\alpha(z) dz= \frac{\lambda^\alpha}{\lambda} \longrightarrow 1, \quad {\alpha \uparrow 1},
			\end{align*}
			then, at the limit $\alpha \uparrow 1$, $\overline{\Pi}^\alpha$ behaves like a delta distribution. This leads us to
			\begin{align*}
				-\int_{x}^{\infty} \frac{d}{dz} f(e,z) \overline{\Pi}^\alpha(z-x) dz \longrightarrow - f^\prime(x), \quad \alpha \uparrow 1.
			\end{align*}
			We can apply the same argument to the vertex $v=(e,0)$ and we have, for $\alpha \uparrow 1$,
			\begin{align*}
				\mathbf{D}^\Phi_{v-} f(v) \longrightarrow - f^\prime_e(v),
			\end{align*}
			where
			\begin{align*}
				f^\prime_e(v)=\lim_{\xi \to 0} \frac{d}{d\xi} f(e,\xi).
			\end{align*}
			Hence, the boundary conditions in Theorem \ref{thm:Xpallino} correspond to the boundary condition of a Walsh Brownian motion when $\alpha$ tends to $1$ in the case of an $\alpha$-stable subordinator. We anticipate that this result holds for the paths of the process as well. The process $B^\bullet$ defined in \eqref{process} approaches a reflecting Brownian motion as $\alpha$ tends to 1 due to $H_t^\alpha \to t$. Consequently, we have $\mathcal{X}^\bullet= (\Theta, B^+)$, where in Definition \ref{def:Xpallino}, $B^\bullet$ is replaced by $B^+$, resulting in a Walsh Brownian motion.
		\end{remark}
		We observed that, also in the case of the star graph, Marchaud (type) derivatives lead us to a process that jumps as soon as it reaches the vertex. The difference compared to the positive half-line is that now the jump can land on other edges based on the $\Theta$ law.
		
		Presently, we want to check if the slowdowns caused by Caputo-D\v{z}rba\v{s}jan (type) derivatives still hold. The non-local sticky Brownian motion on the star graph has been introduced in \cite[Theorem 22]{BonDov}, where at the vertex there are the conditions
		\begin{align*}
			\displaystyle
			c \,\mathfrak{D}_t^\Phi u(t,v)=b \sum_{e \in \mathcal{E}} \rho_e u^\prime_e(t,0),
		\end{align*}
		with  $b+c=1$ and $\mathfrak{D}_t^\Phi$ the Caputo-D\v{z}rba\v{s}jan (type) derivative introduced in \eqref{Caputo}. The aim is to extend this result in the presence of jumps, just as we did for the positive half-line. To achieve this, we introduce the conditions
		\begin{align}
			\varrho, \dot{\varrho} \in C(0,\infty) \text{ and }  \dot{\varrho}(s) \, \overline{\Pi}^\Psi(t-s) \in L^1(0,t),\, t>s>0
			\label{condDL}
		\end{align}
		for a given function $\varrho$ and the space
		\begin{align*}
			D_L:=\left\{\varphi \in C((0,\infty) \times \mathcal{G}) \text{ with } \varrho=\varphi\big\vert_{\mathsf{x}=v} \text{ s.t. $\eqref{condDL}$ are satisfied}   \right\}.
		\end{align*}
		We also define the random times
		\begin{align*}
			\mathcal{S}_t=t+ H^\Psi \circ \left(\eta L^\Phi \circ \ell_t\right)
		\end{align*}
		where $\ell$ is the local time at the vertex of a Walsh Brownian motion on $\mathcal{G}$ with weights $\{\rho_e\}$, and $\mathcal{S}^{-1}$ is the right inverse of $\mathcal{S}$.
		
		We know that $\ell$ coincides with the local time at zero $\gamma$ of the reflected Brownian motion $B^+$ (see \cite[Corollary 8]{BonDov}).
		\begin{theorem}
			\label{thm:Xpallinodelay}
			The probabilistic representation of the solution $u \in D_H \cap D_L$ of the problem
			\begin{align*}
				\begin{cases}
					\dot{u}(t,\mathsf{x})= \Delta u(t,\mathsf{x}) \quad &t>0, \mathsf{x} \in \mathcal{G}\setminus \{v\}, \\
					\eta \mathfrak{D}_t^\Psi u(t,v)+ \sum_{e \in \mathcal{E}} \rho_e\mathbf{D}^\Phi_{x-} u_e(t,x)\big\vert_{x=0}=0 \quad    &t>0,\\
					u(0,\mathsf{x})=f(\mathsf{x}) \quad &\mathsf{x} \in \mathcal{G},
				\end{cases}
			\end{align*}
			with $f \in C(\mathcal{G}) \cap L^\infty(\mathcal{G})$, $\eta \geq 0$  and $\mathsf{x}=(l,x)$ is
			\begin{align*}
				u(t,\mathsf{x})=\mathbf{E}_\mathsf{x}[f(\mathcal{X}^\bullet \circ \mathcal{S}^{-1}_t)].
			\end{align*}
		\end{theorem}
		\begin{proof}
			The proof of this theorem follows that of Theorem \ref{fracdyntm}, where the key observation is that once the edge is specified, the dynamics on $(0,\infty)$ remain the same. Let us examine the principal steps.
			
			We write the $\lambda-$potential, for $\lambda >0$,
			\begin{align*}
				\mathcal{R}_\lambda f(\mathsf{x})=\int_{0}^{\infty} e^{-\lambda t} u(t,\mathsf{x}) dt =\int_{0}^{\infty} e^{-\lambda t} \mathbf{E}_\mathsf{x}[f(\mathcal{X}^\bullet \circ \mathcal{S}^{-1}_t)] dt.
			\end{align*}
			Since $\mathsf{x}=(l,x)$, before hitting the vertex the motion is a killed Brownian motion on the edge $e$. Then, with the same arguments that led us to write \eqref{upsilonSol},  and once a vertex is reached, we choose the edge as in \eqref{resolvent2}, we have, for $\lambda>0$,
			\begin{align}
				\label{resolvent-graph}
				\mathcal{R}_\lambda f(\mathsf{x})=\mathcal{R}_\lambda^D f(\mathsf{x})+ e^{-x \sqrt{\lambda}}  \sum_{e \in \mathcal{E}} \rho_e \mathbf{E}_{(e,0)}  \left[\int_0^\infty e^{-\lambda t} f(e, B^\bullet\circ {S}^{-1}_t) dt\right],
			\end{align}
			where ${S}^{-1}$ is the same random time introduced in Theorem \ref{fracdyntm}, this is possible since the local time $\ell$ at the vertex coincides with $\gamma$.

			As previously noted, the heat equation holds, and we only need to verify the conditions at the vertex. 
			As mentioned earlier, once the edge is specified, the $\lambda-$potential to be computed remains that of the process $B^\bullet\circ {S}^{-1}$ on $(0,\infty)$, which is already determined in \eqref{I1piuI2} through the sum of $I_1^\Psi+ I_2^\Psi$. Then, we have
			\begin{align*}
				&\mathbf{E}_{(e,0)}  \left[\int_0^\infty e^{-\lambda t} f(e, B^\bullet\circ {S}^{-1}_t) dt\right]\\&=\frac{1}{\eta \Psi(\lambda) +\Phi(\sqrt \lambda)} \left( \int_0^\infty R_\lambda^D f(e,y) \, \Pi^\Phi(dy) + \frac{\eta}{\lambda} \Psi(\lambda)f(e,0) \right)
			\end{align*}
			and the $\lambda-$potential at the vertex, by using \eqref{resolvent-graph}, is
			\begin{align*}
				\mathcal{R}_\lambda f(v)= \sum_{e\in \mathcal{E}} \rho_e \left(\frac{1}{\eta \Psi(\lambda) +\Phi(\sqrt \lambda)} \left( \int_0^\infty R_\lambda^D f(e,y) \, \Pi^\Phi(dy) + \frac{\eta}{\lambda} \Psi(\lambda)f(e,0) \right)\right).
			\end{align*}
			We now verify the Laplace transforms for the conditions at the vertex. For the Caputo-D\v{z}rba\v{s}jan (type) derivative, as in \eqref{timeBondCond}, we get
			\begin{align}
				\label{caputo-graphs-bc}
				\int_0^\infty e^{-\lambda t} \eta\mathfrak{D}_t^\Psi u(t,v) dt = \eta \Psi(\lambda) \mathcal{R}_\lambda f(v) - \eta \frac{\Psi(\lambda)}{\lambda} f(v), \quad \lambda>0,
			\end{align}
			where we emphasize that, in exploring solutions within space $D_H$, the function at the vertex is continuous. For the Marchaud (type) derivative, from \eqref{marchaud-graphs}, we have
			\begin{align}
				\label{marchaud-graphs-BC}
				\sum_{e \in \mathcal{E}} \rho_e \mathbf{D}_{x-}^\Phi \mathcal{R}_\lambda f(e,x)\big\vert_{x=0} =-\sum_{e \in \mathcal{E}} \rho_e \int_0^\infty \mathcal{R}_\lambda^D f(e,y)\, \Pi^\Phi(dy) + \Phi(\sqrt{\lambda})\, \mathcal{R}_\lambda f(v).
			\end{align}
			By combining \eqref{caputo-graphs-bc} and \eqref{marchaud-graphs-BC}, the vertex conditions applied to the $\lambda-$potential become
			\begin{align*}
				\eta \Psi(\lambda) \mathcal{R}_\lambda f(v) - \eta \frac{\Psi(\lambda)}{\lambda} f(v) - \sum_{e \in \mathcal{E}} \rho_e \int_0^\infty \mathcal{R}_\lambda^D f(e,y)\, \Pi^\Phi(dy) + \Phi(\sqrt{\lambda})\, \mathcal{R}_\lambda f(v)=0
			\end{align*}
			which are satisfied by
			\begin{align*}
				\mathcal{R}_\lambda f(v)= \sum_{e\in \mathcal{E}} \rho_e \left(\frac{1}{\eta \Psi(\lambda) +\Phi(\sqrt \lambda)} \left( \int_0^\infty R_\lambda^D f(e,y) \, \Pi^\Phi(dy) + \frac{\eta}{\lambda} \Psi(\lambda)f(e,0) \right)\right),
			\end{align*}
			that is the $\lambda-$potential at the vertex of the process $\mathcal{X}^\bullet \circ \mathcal{S}^{-1}$. This concludes that $u(t,\mathsf{x})=\mathbf{E}_\mathsf{x}[f(\mathcal{X}^\bullet \circ \mathcal{S}^{-1}_t)]$ is a solution to the problem under consideration.
		\end{proof}

		In case \( \vert \mathcal{E} \vert = 1 \), we have a description of the paths in terms of \( B^\bullet \circ S^{-1} \) on the half-line, as shown in Figure \ref{StickyBpallino}. We now turn to the general scenario involving the graph \( \mathcal{G} \). According to Definition \ref{def:Xpallino}, the process \( X^\bullet \) is a reflecting Brownian motion on \( \mathcal{G} \) which jumps away from the vertex according with the jump of $B^\bullet$ and the angular process $\Theta$. Once near the vertex, the probability of choosing a specific edge is determined by the weights \( \{\rho_e= \mathbf{P}(\Theta(t) = e)\}_{e \in \mathcal{E}} \). The jumps of $X^\bullet$ well accords with the jumps of \( B^\bullet \). This dynamic is prescribed according with Theorem \ref{thm:Xpallino}, that is case \( \eta = 0 \). We underline that \( \mathcal{T} \), the time at which the process \( \mathcal{X}^\bullet \) hits the vertex, plays the same role as \( \tau_0 \) in the case of the process \( B^\bullet \). For \( \eta > 0 \), as presented in Theorem \ref{thm:Xpallinodelay}, the time change \( \mathcal{S}^{-1} \) introduces the same type of delay as in \( B^\bullet \circ S^{-1} \). Specifically, the process \( \mathcal{X}^\bullet \circ \mathcal{S}^{-1} \) behaves like a Brownian motion on the graph that, upon reaching a vertex, jumps to a new point in \( \mathcal{G} \) (that is, some edge of $\mathcal{G}$), it stops there for a random amount of time according with \( L^\Psi \), the inverse of \( H^\Psi \). Then, it starts afresh.\\
		
		According with Theorem \ref{thm:HTstar} and the discussion around the formula \eqref{PHIderivative}, we get that, for every $i$,  
		\begin{align}
			\mathbf{E}[e^\star_i] = \int_0^\infty e^{-s/\eta} \mathbf{P}(L^\Psi_t \in ds)\, dt = \int_0^\infty \mathbf{P}(L^\Psi_t < \chi)\, dt
		\end{align}
		where $\chi$ is an exponential r.v. with $\mathbf{E}[\chi] = \eta$. That is, for every (visit) $i$, 
		\begin{align}
			\mathbf{E}[e^\star_i] = \int_0^\infty \mathbf{P}(t < H^\Psi \circ \chi)\, dt = \Psi^\prime(0)\, \eta.
			\label{meanHTgraph}
		\end{align}
		We say that the vertex of $\mathcal{G}$ is trap for $\mathcal{X}^\bullet \circ \mathcal{S}^{-1}_t$ if $\Psi^\prime (0)=\infty$. That is, the process spend an infinite average of time on the star vertex.

		\section{Applications}
		
		\subsection{Stochastic resetting}
		In Section \ref{sec:resetting}, we have highlighted a potential connection between NLBVPs and stochastic resetting, a concept widely applied in statistical mechanics, particularly in the context of search problems. We focus on a Brownian particle that intermittently returns to a fixed point (zero in our case), as introduced in \cite{resetting1, resetting2}. The resets occur at random intervals governed by a Poisson process, meaning that the time between successive resets follows an exponential distribution. In this framework, there is a constant probability \(r\) of resetting at any given moment, independent of the time elapsed since the previous reset.
		
		Stochastic resetting proves especially advantageous in search problems due to a key feature: while standard Brownian motion leads to an infinite mean first passage time (MFPT), the introduction of resetting ensures that the MFPT becomes finite \cite{resetting1}.
		
		The process we describe through NLBVPs is closely related to a Brownian motion that resets to zero and subsequently reflects on the positive half-line. This connection with NLBVPs allows us to explore search problems by analyzing repelling boundaries, which induce jumps in the process. These boundaries provide deeper insights into the system's dynamics, offering a richer structure for understanding how resetting and boundary conditions influence behavior.

		The NLBVPs together with the time reversal and the translation of the associated motions give a very large family of processes describing stochastic resetting and stochastic delay. 
		In particular, in \cite{ColDovPag}, distributional equalities are proved both when the initial law is the stationary distribution and also by splitting sample paths at the jump times and concatenating the resulting segments using the Markov property; this pathwise decomposition yields the same distribution for the appropriately constructed time-reversed process.
		
		\subsection{Traffic flow models and delayed agents}
		\label{subsec:traffic-delay}
		
		Sticky diffusions on networks provide a flexible probabilistic framework to model particles (or agents) that spend a finite time at vertices before continuing their motion. 
		Such models have been applied not only to classical road traffic, but also to other transport systems where short-term holding or delays occur at network nodes. 
		In particular, sticky diffusions have been used in the context of Mean Field Games (MFGs) to represent agents that spend a finite (random) time at a vertex before moving on; see \cite{berry2025stationary} for a study of stationary MFGs with sticky transition conditions. The probabilistic setting, see  \cite{berry2024sticky}, assumes exponential holding times at vertices; under this hypothesis the resulting processes are Markovian. 
		However, the modelling framework admits a natural generalisation to non-exponential holding times. 
		In that case the boundary condition at a vertex becomes time-nonlocal (and can be expressed via non-local or fractional time-derivatives at the vertex), which captures memory effects and heavy-tailed waiting-time distributions.
		
		More concretely, by choosing the jump rule according to the edge length one can pass from metric-graph models (particles moving continuously along edges) to models in which particles effectively move between vertices. 
		The time-nonlocal boundary condition then implies random holding times that naturally describe many practical situations: for road traffic, the holding time on an edge may model a broken traffic light or a temporary blockage; for data traffic, it can model a server failure or queueing delay between servers. 
		These examples illustrate how sticky and time-nonlocal boundary mechanisms offer a unified way to incorporate local delays and long-memory effects in network transport models.
		
		

		\subsection{Financial Applications}
		Figure \ref{StickyBpallino} gives an example of possible applications. We may consider motions different from Brownian motions, then our model may entail different properties. Concerning the applications in case we have no jumps (only time non-local operator on the boundary), an example is given in \cite{BKLT-financeSticky} where the authors deal with a model of sticky expectations in which investors update their beliefs too slowly. Our results in this regards give a control for the slow update depending on $\alpha \in (0,1)$. A further example is given by the bank interest rates (or the Vasicek model for instance)  which are termed sticky if they react slowly to changes in the corresponding market rates or in the policy rate (\cite{GNSS-finance-Sticky}). Moreover, in case we have only jumps (only space non-local operator on the boundary) can be given in terms of structural breaks. The process jumps according with the jumps of $B^\bullet_t$ (as in Figure \ref{Bpallino}).

		\subsection{Delayed and rushed reflection}
		
		In the same spirit of the results in \cite{capitanelli2020delayed} we are able to give a classification for the reflection near the boundary. Consider \eqref{meanHoldingTime} and write
		\begin{align*}
			c := \lim_{\lambda \downarrow 0} \frac{\Psi(\lambda)}{\lambda}
		\end{align*}
		assuming the limit exists ($c=\infty$ is a clear case, $\mathbf{E}[H^\Psi_t] = \infty$). The process $\bar{X}$ may have (in terms of mean holding time on the boundary):
		\begin{itemize}  
			\item[B1)] delayed boundary behavior if $c<1$;
			\item[B2)] rushed boundary behavior if $c>1$;
			\item[B3)] base (reflected) boundary behavior if $c=1$.
		\end{itemize}
		A special case of subordinator including either delayed or rushed effect is the gamma subordinator for which, for $a,b>0$,
		\begin{align*}
			\Psi(\lambda) = a \ln (1+ \lambda/b) \quad \textrm{and} \quad \lim_{\lambda \downarrow 0} \frac{\Psi(\lambda)}{\lambda} = \frac{a}{b} 
		\end{align*}
		and $c<1$, $c>1$ or $c=1$ depending on the ratio $a/b$. On the other hand, by considering $\Psi=Id$, we have that $B^\bullet$ jumps away from the boundary according with the last jump of $H^\Phi$. Due to nature of the subordinator, the process $B^\bullet$ never hits the boundary, it reflects instantaneously before hitting the boundary. In this case, the time that the process spends on the boundary can be only given (eventually) by the starting point as a point of the boundary, case in which the process is pushed immediately away once again for the nature of the subordinator. We recall that our focus is only on strictly increasing subordinators. 
		
		Our analysis introduces the following further characterization: 
		\begin{itemize}  
			\item[R1)] delayed/rushed (continuous) slow reflection;
			\item[R2)] jumping (discontinuous) reflection;
			\item[R3)] (continuous) reflection.
		\end{itemize}
		The delayed reflection can be related with a reflection on the boundary (like an insulating boundary for example) whereas the rushed reflection seems to be related with a reflection near the boundary (like an inflating boundary for example). These characterizations are therefore given according with the comparison between the total (mean) amounts of time the processes $B^+$ and $\bar{X}$ would spend on the boundary. By base boundary behavior we mean the behavior of the reflecting process $B^+$. We observe that instantaneous reflection means that the time the process spends on the boundary has zero Lebesgue measure. The reflection R1 is obtained via time change (by considering $S^{-1}_t$) and depends on the symbol $\Psi$ of the boundary condition. R1 gives rise to the behaviors B1 and B2, notice that in this case we may consider a boundary diffusion (in the sense of Dirichlet-Neumann generator) which can be delayed or rushed with respect to the base boundary diffusion (according with the definition given in \cite{capitanelli2020delayed}). R2 turns out to be associated with the symbol $\Phi$ of the boundary condition. Concerning R3 we easily make the connection with B3, that is the case of Neumann boundary condition.\\
		
		In case $E=[0, \infty)$ we are considering the trapping point $\{0\}$. Despite of the apparently simple situation, it basically represents a very general setting.

		%
		%
		%
		%
		
		\section{Proof of Theorem \ref{fracdyntm}}
		\label{proof1}
		By exploiting the Markov property of the Brownian motion, as a solution to the heat equation, the function $\upsilon$ can be written as
		\begin{align}
			\label{abWITHsol}
			\upsilon(t,x) = a(t,x) + \int_0^t b(t-s, x) \upsilon(s,0)ds
		\end{align}
		with Laplace transform of both side
		\begin{align*}
			\tilde{\upsilon}(\lambda, x) = \tilde{a}(\lambda, x) + \tilde{b}(\lambda, x)\, \tilde{\upsilon}(\lambda, 0), \quad \lambda>0
		\end{align*}
		for some suitable (sufficiently regular and integrable) functions $a,b$. Indeed, simple arguments say that $\tilde{a}(\lambda, 0)=0$ and $\tilde{b}(\lambda, 0)=1$. Moreover, it must be $\upsilon(0,x) = a(0, x)$ and $\tilde{\upsilon}^{\prime \prime} = \lambda \tilde{\upsilon} - f$. In particular, it turns out that $a,b$ are given by
		\begin{align}
			\label{upsilonSol}
			\upsilon(t,x)=Q_t^D f(x)+\int_0^t \frac{x}{\tau} g(\tau,x) \upsilon(t-\tau,0) d\tau
		\end{align}
		where
		\begin{align*}
			Q_t^D f(x)=\int_0^\infty (g(t,x-y)-g(t,x+y))f(y) dy
		\end{align*}
		is the Dirichlet semigroup and $g(t, z) = e^{-z^2/ 4t}/\sqrt{4 \pi t}$  for which we recall the well-known transforms
		\begin{align}
			\label{gLap1}
			\int_0^\infty e^{-\lambda t} g(t,x) dt = \frac{1}{2}\frac{e^{-x \sqrt \lambda}}{\sqrt \lambda}, \quad \lambda >0
		\end{align}
		and
		\begin{align}
			\label{gLap2}
			\int_0^\infty e^{-\lambda t} \frac{x}{t} g(t,x) dt = e^{-x \sqrt \lambda}, \quad \lambda>0.
		\end{align}
		In particular, if $a$ is the Dirichlet semigroup, then $b$ must be obtained from the unique (continuous and bounded) solution to $\tilde{b}^{\prime \prime} = \lambda \tilde{b}$ with $\tilde{b}(\lambda, 0)=1$. Thus, $\upsilon(t, 0)$ in \eqref{upsilonSol} completely identifies the boundary condition and vice versa.

		Denote by $\tau^\bullet$ the first time $B^\bullet$ hits the boundary point $0$ and observe that $\mathbf{P}_x(B^\bullet_{\tau^\bullet-}=0)=1$, $\forall\, x>0$. From the probability viewpoint, first we notice that 
		\begin{align*}
			R_1(x) 
			= & \mathbf{E}_x \left[ \int_0^\infty e^{-\lambda t} f(B^\bullet_t) dt \right]\\ 
			= & \mathbf{E}_x \left[ \int_0^{\tau^\bullet} e^{-\lambda t} f(B^\bullet_t) dt \right] + \mathbf{E}_x \left[ \int_{\tau^\bullet}^\infty e^{-\lambda t} f(B^\bullet_t) dt \right]
		\end{align*}
		and
		\begin{align*}
			R_2(x) 
			= & \mathbf{E}_x \left[ \int_0^\infty e^{-\lambda t} f(B^\bullet \circ S^{-1}_t) dt \right]\\ 
			= & \mathbf{E}_x \left[ \int_0^{\tau^\bullet} e^{-\lambda t} f(B^\bullet_t) dt \right] + \mathbf{E}_x \left[ \int_{\tau^\bullet}^\infty e^{-\lambda t} f(B^\bullet \circ S^{-1}_t) dt \right]
		\end{align*}
		respectively give the resolvent operator for the Markov process $B^\bullet$ and the $\lambda$-potential for $B^\bullet \circ S^{-1}_t$. This well-agrees with formula \eqref{abWITHsol}. Both $R_1$ and $R_2$ involve a Dirichlet semigroup for $B^\bullet$ whereas the second part of $R_2$ writes
		\begin{align*}
			\mathbf{E}_x \left[ \int_{\tau^\bullet}^\infty e^{-\lambda t} f(B^\bullet \circ S^{-1}_t) dt \right] 
			= & \mathbf{E}_x \left[ \int_{S \circ \tau^\bullet}^\infty e^{-\lambda S_t} f(B^\bullet_t) dS_t \right]\\
			= & - \frac{1}{\lambda} \mathbf{E}_x \left[ \int_{S \circ \tau^\bullet}^\infty  f(B^\bullet_t) de^{-\lambda S_t} \right].
		\end{align*}
		We can invoke the Markov property of $B^\bullet$ in order to confirm the analytic discussion on \eqref{abWITHsol}. We will provide conclusive arguments on Section \ref{sec:HTandRT}. Now we focus on the probabilistic representation of the solution $\upsilon$.\\
		
		Let us write
		\begin{align}
			\label{upsilonASS}
			R_\lambda f(x) = \int_0^\infty e^{-\lambda t} \upsilon(t,x) dt, \quad \lambda >0.
		\end{align}		
		Our aim is to prove that
		\begin{align}
			\label{potentialASS}
			R_\lambda f(x)=\mathbf{E}_x\left[\int_0^\infty e^{-\lambda t} f(B^\bullet \circ S_t^{-1}) dt\right], \quad  \lambda >0
		\end{align}
		is the $\lambda$-potential associated with the process $B^\bullet \circ S^{-1}$.
		
		Under the representation \eqref{upsilonSol}, formula \eqref{upsilonASS} takes the form
		\begin{align*}
			R_\lambda f(x) = R^D_\lambda f(x) + \bar{R}_\lambda f(x)
		\end{align*}
		where $R^D_\lambda f$ and $\bar{R}_\lambda f$ are defined below. From \eqref{gLap1} and \eqref{gLap2} we respectively write
		\begin{align}
			\label{Rtemp}
			R^D_\lambda f(x) 
			& =\int_0^\infty e^{-\lambda t } Q_t^D f(x) dt\notag \\	
			&=\frac{1}{2}\int_0^\infty \left(\frac{e^{-\vert x-y \vert \sqrt \lambda}}{\sqrt \lambda}- \frac{e^{-( x+y)  \sqrt \lambda}}{\sqrt \lambda}\right) f(y) dy\\
			&=\frac{1}{2} \int_0^\infty \frac{e^{( x-y)  \sqrt \lambda}}{\sqrt \lambda} f(y) dy - \frac{1}{2} \int_0^\infty \frac{e^{-( x+y)  \sqrt \lambda}}{\sqrt \lambda} f(y) dy\ +\notag\\
			&- \frac{1}{2} \int_0^x \left(\frac{e^{( x-y)  \sqrt \lambda}}{\sqrt \lambda} - \frac{e^{-( x-y)  \sqrt \lambda}}{\sqrt \lambda} \right) f(y) dy, \quad \lambda>0 \notag
		\end{align}
		and
		\begin{align*}
			\bar{R}_\lambda f(x) = e^{-x \sqrt{\lambda}} R_\lambda f(0), \quad \lambda>0.
		\end{align*}
		We can easily check that
		\begin{align}
			\label{resolvent_heat}
			\big(R_\lambda f(x) \big)^{\prime \prime}= \lambda R_\lambda f(x) - f(x) = \int_0^\infty e^{-\lambda t} \frac{\partial}{\partial t} \upsilon(t, x)dt, \quad \lambda>0, x >0
		\end{align}
		and $\upsilon$ solves the heat equation. We also notice that $R_\lambda f$ is continuous at $x=0$ and $R_\lambda f(0) = \bar{R}_\lambda f(0)$. The characterization of $R_\lambda f(0)$ and therefore $\upsilon(t,0)$ gives \eqref{upsilonSol}. Thus, we focus on the boundary condition for which, first we observe that
		\begin{align*}
			\int_0^\infty e^{-\lambda t} \mathfrak{D}_t^\Psi \upsilon(t,x)\bigg\vert_{x=0} dt  = \int_0^\infty e^{-\lambda t} \mathfrak{D}_t^\Psi \upsilon(t,x) dt \bigg\vert_{x=0}, \quad \lambda>0
		\end{align*}
		and, thanks to \eqref{LapCaputo}, we have
		\begin{align}
			\label{timeBondCond}
			\int_0^\infty e^{-\lambda t} \eta\mathfrak{D}_t^\Psi \upsilon(t,x) dt \bigg\vert_{x=0}= \eta \Psi(\lambda) R_\lambda f(0) - \eta \frac{\Psi(\lambda)}{\lambda} f(0), \quad \lambda>0.
		\end{align}
		A key point in our proof is the linearity of $\mathbf{D}_{x-}^\Phi$ for which 
		\begin{align*}
			\mathbf{D}_{x-}^\Phi R_\lambda f(x) \big\vert_{x=0} = \big(\mathbf{D}_{x-}^\Phi R_\lambda^D f(x) + \mathbf{D}_{x-}^\Phi \bar{R}_\lambda f(x)\big)\big\vert_{x=0}
		\end{align*}
		as we can check from the definition \eqref{Marchaudleft}. We use the fact that
		\begin{align*}
			\lim_{x\downarrow 0}  \int_0^\infty \left( R^D_\lambda f(x) - R^D_\lambda f(x+y) \right) \Pi^\Phi(dy) = - \int_0^\infty R^D_\lambda f(y)\, \Pi^\Phi(dy)
		\end{align*}
		and
		\begin{align*}
			&\lim_{x\downarrow 0}  \int_0^\infty \left( \bar{R}_\lambda f(x) - \bar{R}_\lambda f(x+y) \right) \Pi^\Phi(dy) \\
			&=  \lim_{x\downarrow 0}  \int_0^\infty \left( e^{-x \sqrt{\lambda}} - e^{-(x+y) \sqrt{\lambda}} \right) \Pi^\Phi(dy)\, R_\lambda f(0)\\
			&= \Phi(\sqrt{\lambda})\, R_\lambda f(0).
		\end{align*}
		Notice that, for $\vert f(x) \vert \leq K$, 
		\begin{align*}
			\int_0^\infty R^D_\lambda f(y)\, \Pi^\Phi(dy)\leq \frac{K}{\lambda}  \int_0^\infty (1-e^{- \sqrt{\lambda} y}) \, \Pi^\Phi(dy)= \frac{K}{\lambda} \Phi(\sqrt{\lambda}), \quad \lambda >0.
		\end{align*}
		We therefore obtain
		\begin{align*}
			- \mathbf{D}_{x-}^\Phi R_\lambda f(x) \big\vert_{x=0} = \int_0^\infty R^D_\lambda f(y)\, \Pi^\Phi(dy) - \Phi(\sqrt{\lambda})\, R_\lambda f(0).
		\end{align*}
		This together with \eqref{timeBondCond} leads to
		\begin{align}
			\label{BCresolvent2}
			R_\lambda f(0) =\frac{1}{\eta \Psi(\lambda) +\Phi(\sqrt \lambda)} \left( \frac{\eta}{\lambda} \Psi(\lambda)f(0)+ \int_0^\infty R_\lambda^D f(y) \, \Pi^\Phi(dy) \right), \quad \lambda>0
		\end{align}
		which has been obtained under the representation \eqref{upsilonSol}.\\
		
		Now we study the probabilistic representation of $R_\lambda f(0)$. Assume that \eqref{potentialASS} holds true. Recall that $S_t = t + H^\Psi \circ (\eta L^\Phi_{\gamma_t})$ is strictly increasing. The inverse $S^{-1}$ may have plateaus. Observe that $R_\lambda f$ can be written in terms of the path integral
		\begin{align*}
			R_\lambda f(0) 
			& = \mathbf{E}_0\left[\int_0^\infty e^{-\lambda t} f(B^\bullet \circ S_t^{-1}) dt\right]\\
			& = \mathbf{E}_0\left[\int_0^\infty e^{-\lambda S_t} f(B^\bullet_t) dS_t\right]\\ 
			& = I_1^\Psi+I_2^\Psi
		\end{align*}
		where 
		\begin{align*}
			I_1^\Psi = \mathbf{E}_0\left[\int_0^\infty e^{-\lambda t}e^{- \lambda H^\Psi\eta  L^\Phi_{\gamma_t}} f(B^\bullet_t) dt\right], \quad \lambda>0
		\end{align*}
		and
		\begin{align*}
			I_2^\Psi = f(0) \mathbf{E}_0\left[\int_0^\infty e^{-\lambda t} e^{-\lambda H^\Psi\eta  L^\Phi_{\gamma_t}} dH^\Psi \eta L^\Phi_{\gamma_t}\right], \quad \lambda>0.
		\end{align*}
		The computation of $I_1^\Psi$ follows from the computation of $e_1$ in \cite[page 215]{ito1963brownian}, recalling that 
		\begin{align*}
			\mathbf{E}_0[e^{-\lambda H^\Psi_{\eta t}}]=e^{-\eta t \Psi(\lambda)}, \quad \lambda>0.
		\end{align*}
		Indeed, $H^\Phi$ is independent from $B^\bullet$ (and therefore, from $L^\Phi$) and we get
		\begin{align*}
			I^\Psi_1 
			= \mathbf{E}_0\left[\int_0^\infty e^{-\lambda t} e^{- \eta \Psi(\lambda)  L^\Phi_{\gamma_t}} f(B^\bullet_t) dt\right].
		\end{align*}
		From  \cite[Formula 6, Section 15]{ito1963brownian}, see the Appendix, we know 
		\begin{align*}
			\int_0^\infty e^{-\lambda t} \mathbf{E}_0[f(B^\bullet_t)]= \frac{\int_0^\infty {R}_\lambda^D f(y) \Pi^\Phi(dy)}{\Phi(\sqrt{\lambda})},
		\end{align*}
		where \({R}_\lambda^D f\) is the resolvent of the Brownian motion on the half-line with the Dirichlet boundary condition at zero. Then, we conclude that
		\begin{align*}
			I_1^\Psi=\frac{1}{\eta \Psi(\lambda) +\Phi(\sqrt \lambda)} \int_0^\infty R_\lambda^D f(y) \, \Pi^\Phi(dy), \quad \lambda>0.
		\end{align*}
		For the path integral $I_2^\Psi$, we get that
		\begin{align*}
			I_2^\Psi&= f(0) \mathbf{E}_0\left[\int_0^\infty e^{-\lambda t} e^{-\lambda H^\Psi\eta  L^\Phi_{\gamma_t}} d(H^\Psi \circ (\eta L^\Phi_{\gamma_t}))\right]\\
			&= f(0)  \mathbf{E}_0\left[\int_0^\infty e^{-\gamma^{-1} H^\Phi_t} e^{-\lambda H^\Psi_{\eta t}} dH^\Psi_{ \eta t}\right]\\
			&=\frac{f(0)}{\lambda} \mathbf{E}_0\left[-1-\int_0^\infty e^{-\lambda H^\Psi_{\eta t} } d\left(e^{-\gamma^{-1} H^\Phi_t}\right)\right]\\
			&=\frac{f(0)}{\lambda} \mathbf{E}_0\left[-1-\int_0^\infty e^{-\eta t \Psi(\lambda) } d\left(e^{-\gamma^{-1} H^\Phi_t}\right)\right]\\
			&=\frac{f(0)}{\lambda} \mathbf{E}_0\left[-1-\left(-1-\eta \Psi(\lambda)\int_0^\infty e^{-\eta t \Psi(\lambda) } e^{-\lambda \gamma^{-1}{H^\Phi_t}} dt \right)\right]\\
			&=\frac{\eta}{\lambda} \Psi(\lambda) f(0) \mathbf{E}_0\left[\int_0^\infty e^{-\eta t \Psi(\lambda)}e^{-\lambda \gamma^{-1} H^\Phi_t} dt\right]\\
			&=\frac{\eta}{\lambda} \Psi(\lambda) f(0) \mathbf{E}_0\left[\int_0^\infty e^{-\eta  \Psi(\lambda) L^\Phi_{\gamma_t}}e^{-\lambda t} d L^\Phi_{\gamma_t}\right], \quad \lambda>0
		\end{align*}
		where
		\begin{align*}
			& \mathbf{E}_0\left[\int_0^\infty e^{-\lambda t} e^{-\eta  \Psi(\lambda) L^\Phi_{\gamma_t}} d L^\Phi_{\gamma_t}\right]\\
			&=-\frac{1}{\eta \Psi(\lambda)} \mathbf{E}_0\left[\int_0^\infty e^{-\lambda t} de^{-\eta  \Psi(\lambda) L^\Phi_{\gamma_t}}\right]\\
			&=-\frac{1}{\eta \Psi(\lambda)}\mathbf{E}_0\left[-1+\lambda  \int_0^\infty e^{-\lambda t } e^{-\eta \Psi(\lambda) L^\Phi_{\gamma_t}} dt\right]\\
			&=\frac{1}{\eta \Psi(\lambda)}-\frac{\lambda}{\eta \Psi(\lambda)} \mathbf{E}_0\left[\int_0^\infty e^{-\lambda t} e^{- \eta \Psi(\lambda) L^\Phi_{\gamma_t}} dt\right]\\
			&=[\textrm{from \eqref{gLap1}}] = \frac{1}{\eta \Psi(\lambda)}-\frac{\lambda}{\eta \Psi(\lambda)} \mathbf{E}_0\left[\int_0^\infty \frac{\sqrt \lambda}{\lambda }e^{-\sqrt \lambda w} e^{- \eta \Psi(\lambda) L^\Phi_w} dw\right]\\
			&= [\textrm{from \eqref{LapL}}] = \frac{1}{\eta \Psi(\lambda)}-\frac{\lambda}{\eta \Psi(\lambda)} \int_0^\infty \frac{\sqrt \lambda}{\lambda } e^{-\eta \Psi(\lambda) z} \frac{\Phi(\sqrt \lambda)}{\sqrt \lambda} e^{-z \Phi(\sqrt \lambda)} dz\\
			&=\frac{1}{\eta \Psi(\lambda)}-\frac{\Phi(\sqrt \lambda)}{\eta \Psi(\lambda)}\frac{1}{\eta \Psi(\lambda)+ \Phi(\sqrt \lambda)}\\
			&=\frac{1}{\eta \Psi(\lambda)+ \Phi(\sqrt \lambda)}.
		\end{align*}
		We arrive at
		\begin{align*}
			I_2^\Psi =  \frac{1}{\eta \Psi(\lambda)+ \Phi(\sqrt \lambda)} \frac{\eta}{\lambda} \Psi(\lambda) f(0), \quad \lambda>0.
		\end{align*}
		By collecting all pieces together, we obtain
		\begin{align}
			I_1^\Psi+I_2^\Psi
			&=\frac{1}{\eta \Psi(\lambda) +\Phi(\sqrt \lambda)} \left( \int_0^\infty R_\lambda^D f(y) \, \Pi^\Phi(dy) + \frac{\eta}{\lambda} \Psi(\lambda)f(0) \right)
			\label{I1piuI2}
		\end{align}
		which coincides with \eqref{BCresolvent2}. Thus, we conclude that \eqref{potentialASS} holds true. In particular, $\upsilon$ is a continuous function with Laplace transform $R_\lambda f$, then $\upsilon$ is unique. There exists only one continuous inverse to $R_\lambda f$.
		
		We also stress the fact that we deal with positive functions on unbounded domain. Standard arguments guarantees uniqueness of the solution in terms of positivity (Widder's theorem \cite{widder1944positive}) and exponential growth (\cite[section 2.3.3]{evans2010partial}).
		\section{Proof of Theorem \ref{thm:Xpallino}}
		\label{proof2}
		Let us introduce the resolvent, for $\lambda >0$,
		\begin{align*}
			\mathcal{R}_\lambda f(\mathsf{x})=\int_{0}^{\infty} e^{-\lambda t} u(t,\mathsf{x}) dt =\int_{0}^{\infty} e^{-\lambda t} \mathbf{E}_\mathsf{x}[f(\mathcal{X}^\bullet_t)] dt.
		\end{align*}
		Since $\mathcal{X}^\bullet$ starts afresh with a jump after $\mathcal{T}$, by proceeding as \cite[Formula 2, Section 15]{ito1963brownian}, we obtain
		\begin{align}
			\label{resolvent}
			\mathcal{R}_\lambda f(\mathsf{x})=\mathbf{E}_\mathsf{x} \left[\int_{0}^{\mathcal{T}} e^{-\lambda t} f(\mathcal{X}^\bullet_t) dt \right] + \mathbf{E}_\mathsf{x}\left[e^{-\lambda \mathcal{T}}\right] \mathbf{E}_v\left[\int_{0}^{\infty} e^{-\lambda t} f(\mathcal{X}^\bullet_t) dt  \right],
		\end{align}
		where now the role of the origin on the positive half-line is taken by the vertex $v$. A key observation is that by expressing the process $\mathcal{X}^\bullet$ as $(\Theta, B^\bullet)$, then the first-hitting time of the vertex is equivalent to the first-hitting time of $B^\bullet$ at zero. For the first part, in the same way of \eqref{Rtemp}, since the process $B^\bullet$ before hitting zero is a killed Brownian motion $B^D$, associated to a Dirichlet boundary condition, we get
		\begin{align*}
			\mathbf{E}_\mathsf{x} \left[\int_{0}^{\mathcal{T}} e^{-\lambda t} f(\mathcal{X}^\bullet_t) dt \right]&=\mathbf{E}_{(l,x)}\int_0^{\tau_0} e^{-\lambda t} f(l, B^\bullet_t)dt\\&=\mathbf{E}_{(l,x)}\int_0^\infty e^{-\lambda t} f(l, B^D_t)dt\\
			&=\int_0^\infty e^{-\lambda t} \int_0^\infty (g(t,x-y)-g(t,x+y))f(l,y) dy \,dt\\
			&=:\mathcal{R}_\lambda^D f(l,x),
		\end{align*}
		where $g(t, z) = e^{-z^2/ 4t}/\sqrt{4 \pi t}$. By using \eqref{gLap1}
		we conclude that
		\begin{align}
			\label{resolvent1}
			&\mathcal{R}^D_\lambda f(l,x)=\mathbf{E}_{(l,x)}\int_0^{T} e^{-\lambda t} f(l, B^\bullet_t) dt \notag \\
			&=\frac{1}{2}\int_0^\infty \left(\frac{e^{-\vert x-y \vert \sqrt \lambda}}{\sqrt \lambda}- \frac{e^{-( x+y)  \sqrt \lambda}}{\sqrt \lambda}\right) f(l,y) dy\notag\\
			&=\frac{1}{2} \int_0^\infty \frac{e^{( x-y)  \sqrt \lambda}}{\sqrt \lambda} f(l,y) dy - \frac{1}{2} \int_0^\infty \frac{e^{-( x+y)  \sqrt \lambda}}{\sqrt \lambda} f(l,y) dy\ +\notag\\
			&- \frac{1}{2} \int_0^x \left(\frac{e^{( x-y)  \sqrt \lambda}}{\sqrt \lambda} - \frac{e^{-( x-y)  \sqrt \lambda}}{\sqrt \lambda} \right) f(l,y) dy, \quad \lambda>0.
		\end{align}
		For the second part of \eqref{resolvent}, since $\Theta \sim \mu$, we have
		\begin{align*}
			\mathbf{E}_\mathsf{x}\left[e^{-\lambda \mathcal{T}}\right] \mathbf{E}_v\left[\int_{0}^{\infty} e^{-\lambda t} f(\mathcal{X}^\bullet_t) dt  \right]&=\mathbf{E}_{(l, x)}\left[e^{-\lambda \mathcal{T}}\right] \sum_{e \in \mathcal{E}} \rho_e \, \mathbf{E}_{(e,0)}  \left[\int_0^\infty e^{-\lambda t} f(e, B^\bullet_t) dt\right]\\
			&= \mathbf{E}_{(l, x)}\left[e^{-\lambda \mathcal{T}}\right] \sum_{e \in \mathcal{E}} \rho_e {\mathcal{R}}_\lambda f(e,0).
		\end{align*}
		We note that $\mathbf{E}_{(l,x)}[e^{-\lambda T}]=\mathbf{E}_{x}[e^{-\lambda \tau_0}]= e^{-x \sqrt{\lambda}}$ for all $l \in \mathcal{E}$, since, once the edge is fixed, it remains a simple Brownian excursion. Then, we get
		\begin{align}
			\label{resolvent2}
			\mathbf{E}_\mathsf{x}\left[e^{-\lambda \mathcal{T}}\right] \mathbf{E}_v\left[\int_{0}^{\infty} e^{-\lambda t} f(\mathcal{X}^\bullet_t) dt\right] =  e^{-x \sqrt{\lambda}}  \sum_{e \in \mathcal{E}} \rho_e \mathbf{E}_{(e,0)}  \left[\int_0^\infty e^{-\lambda t} f(e, B^\bullet_t) dt\right].
		\end{align}
		By combining \eqref{resolvent1} and  \eqref{resolvent2}, we see that, for the resolvent \eqref{resolvent}, the following holds
		\begin{align*}
			\lambda \mathcal{R}_\lambda f(l,x) - f(l, x)= \frac{\partial^2}{\partial x^2} \mathcal{R}_\lambda f(l,x),
		\end{align*}
		then the heat equation 
		\begin{align*}
			\begin{cases}
				\dot{u}(t,\mathsf{x})= \Delta u(t,\mathsf{x}) \quad &t>0, \mathsf{x} \in \mathcal{G}\setminus \{v\} \\
				\upsilon(0,\mathsf{x})=f(\mathsf{x}) \quad &\mathsf{x} \in \mathcal{G}
			\end{cases}
		\end{align*}
		is satisfied. Now, let us move on the boundary conditions. When the process is in $\mathsf{x}=v$, the resolvent simply reduces to
		\begin{align*}
			\mathcal{R}_\lambda f(v)= \sum_{e \in \mathcal{E}} \rho_e \mathbf{E}_{(e,0)}  \left[\int_0^\infty e^{-\lambda t} f(e, B^\bullet_t) dt\right].
		\end{align*}
		From  \cite[Formula 6, Section 15]{ito1963brownian}, see the Appendix, we know that on the positive half-line for \(B^\bullet\)
		\begin{align*}
			\int_0^\infty e^{-\lambda t} \mathbf{E}_0[f_e(B^\bullet_t)]= \frac{\int_0^\infty {R}_\lambda^D f_e(l) \Pi^\Phi(dl)}{\Phi(\sqrt{\lambda})},
		\end{align*}
		where \({R}_\lambda^D f\) is the resolvent of the Brownian motion on the half-line with the Dirichlet boundary condition at zero. From the definition of \(\mathcal{X}^\bullet\), we conclude that the resolvent in \(v\) is
		\begin{align*}
			\mathcal{R}_\lambda f(v)&= \sum_{e \in \mathcal{E}} \rho_e \mathbf{E}_{(e,0)}  \left[\int_0^\infty e^{-\lambda t} f(e, B^\bullet_t) dt\right]\\&=\sum_{e \in \mathcal{E}} \rho_e \int_0^\infty e^{-\lambda t} \mathbf{E}_0[f_e(B^\bullet_t)]\\&= \sum_{e \in \mathcal{E}} \rho_e \frac{\int_0^\infty \mathcal{R}_\lambda^D f(e, l) \Pi^\Phi(dl)}{\Phi(\sqrt{\lambda})}.
		\end{align*}
		For more details on the resolvent in the vertex, see also \cite[Section 21.11]{werner2016brownian}.
		
		
		Now we show that the vertex conditions have the same Laplace transform.
		For the Marchaud (type) derivative, fixed an $l \in \mathcal{E}$, we see that
		\begin{align*}
			\sum_{e \in \mathcal{E}}\rho_e \mathbf{D}_{x-}^\Phi \mathcal{R}_\lambda f(e,x)\big\vert_{x=0}  = \left(\sum_{e \in \mathcal{E}} \rho_e \mathbf{D}_{x-}^\Phi \mathcal{R}_\lambda^D f(e,x)\big\vert_{x=0} + \mathbf{D}_{x-}^\Phi e^{-\lambda x}\big\vert_{x=0}\sum_{e \in \mathcal{E}} \rho_e {\mathcal{R}}_\lambda f(e,0)\right)
		\end{align*}
		In particular, for the killed part
		\begin{align*}
			&\sum_{e \in \mathcal{E}} \rho_e \lim_{x\downarrow 0}  \int_0^\infty \left( \mathcal{R}^D_\lambda f(e,x) - \mathcal{R}^D_\lambda f(e,x+y) \right) \Pi^\Phi(dy) \\
			& = \sum_{e \in \mathcal{E}} \rho_e \int_0^\infty \left( \mathcal{R}^D_\lambda f(v) - \mathcal{R}^D_\lambda f(e,y) \right) \Pi^\Phi(dy) \\
			& = - \sum_{e \in \mathcal{E}} \rho_e \int_0^\infty \mathcal{R}^D_\lambda f(e,y)\, \Pi^\Phi(dy)
		\end{align*}
		and
		\begin{align*}
			\lim_{x\downarrow 0}  \int_0^\infty \left( e^{-x \sqrt{\lambda}} - e^{-(x+y) \sqrt{\lambda}} \right) \Pi^\Phi(dy)\, \sum_{e \in \mathcal{E}} \rho_e {\mathcal{R}}_\lambda f(e,0)
			&= \Phi(\sqrt{\lambda})\, \sum_{e \in \mathcal{E}} \rho_e {\mathcal{R}}_\lambda f(e,0).
		\end{align*}
		We therefore obtain that
		\begin{align}
			\label{marchaud-graphs}
			0&=- \sum_{e \in \mathcal{E}} \rho_e \mathbf{D}_{x-}^\Phi \mathcal{R}_\lambda f(e,x)\big\vert_{x=0}\\
			& =\sum_{e \in \mathcal{E}} \rho_e \int_0^\infty \mathcal{R}_\lambda^D f(e,y)\, \Pi^\Phi(dy) - \Phi(\sqrt{\lambda})\, \sum_{e \in \mathcal{E}} \rho_e {\mathcal{R}}_\lambda f(e,0),
		\end{align}
		it is satisfied by
		\begin{align}
			\label{BCresolvent2G}
			\mathcal{R}_\lambda f(v) = \sum_{e \in \mathcal{E}} \rho_e {\mathcal{R}}_\lambda f(e,0)=\sum_{e \in \mathcal{E}} \rho_e \frac{\int_0^\infty \mathcal{R}_\lambda^D f(e, l) \Pi^\Phi(dl)}{\Phi(\sqrt{\lambda})}, \quad \lambda>0.
		\end{align}
		Then, we provide that $u(t,\mathsf{x})=\mathbf{E}_\mathsf{x}[f(\mathcal{X}^\bullet_t)]$ solves the boundary conditions
		\begin{align*}
			\sum_{e \in \mathcal{E}} \rho_e  \mathbf{D}_{x-}^\Phi f(e,x)\big\vert_{x=0}=0.
		\end{align*}
		and since we have already seen that it solves the heat equation, it is the probabilistic solution of \eqref{problem}.

		\section{Appendix}
		
		\label{Sec:Appendix}

		\subsection{Some auxiliary results}
		
		\label{sec:aux}
		
		Observe that $AC(I)$ coincides with $W^{1,1}(I)$ only if $I\subset (0,\infty)$ is bounded. We recall that $AC$ denotes the set of absolutely continuous functions. For the interval $I$, $W^{1,\infty}(I)$ coincides with the space of Lipschitz continuous functions.

		In order to define the Riemann-Liouville (type) derivatives on the positive real line, we first consider a closed interval $\bar{I} \subset (0,\infty)$ and $u \in AC(\bar{I})$. 
		Then we extend the result on $\mathbb{R}^+$ as in \cite[page 79]{kilbas2006theory}. We now introduce the Riemann-Liouville (type) derivatives
		\begin{align}
			\label{R-Lxinf}
			\mathcal{D}_{(x,\infty)}^\Phi u(x):=-\frac{d}{dx} \int_x^\infty u(y) \overline{\Pi}^\Phi(y-x) dy
		\end{align}
		and
		\begin{align}
			\label{R-L0x}
			\mathcal{D}_{(0,x)}^\Phi u(x):=\frac{d}{dx} \int_0^x u(y) \overline{\Pi}^\Phi(x-y) dy
		\end{align}
		respectively defined for function $u$ such that 
		\begin{align*}
			u(\cdot) \overline{\Pi}^\Phi(\cdot-x) \in L^1(x,\infty), \quad \text{and} \quad u(\cdot) \overline{\Pi}^\Phi(x-\cdot) \in L^1(0,x) \quad \forall x.
		\end{align*}
		Let us focus on \eqref{R-Lxinf}. We observe that
		\begin{align*}
			\mathcal{D}_{(x,\infty)}^\Phi u(x)	&=-\frac{d}{dx} \int_x^\infty u(y) \overline{\Pi}^\Phi(y-x) dy\\&=\lim_{b \to \infty} -\frac{d}{dx} \int_x^b u(y) \overline{\Pi}^\Phi(y-x) dy\\
			&=\lim_{b \to \infty} -\frac{d}{dx} \int_0^{b-x} u(z+x) \overline{\Pi}^\Phi(z) dz\\
			&=\lim_{b \to \infty} u(b) \overline{\Pi}^\Phi(b-x) - \int_0^{b-x} u^\prime(z+x) \overline{\Pi}^\Phi(z) dz\\
			&=\lim_{b \to \infty} u(b) \overline{\Pi}^\Phi(b-x) - \int_x^{\infty} u^\prime(y) \overline{\Pi}^\Phi(y-x) dy
		\end{align*}
		where, from the final value theorem for the Laplace transform in \eqref{tailSymb}, we have
		\begin{align*}
			0=\Phi(0)=\lim_{b \to \infty} \overline{\Pi}^\Phi(b).
		\end{align*}
		Assuming that the growth of $u$ is asymptotically bounded, then
		\begin{align}
			\label{R-LCaputo}
			\mathcal{D}_{(x,\infty)}^\Phi u(x)=- \int_x^{\infty} u^\prime(y) \overline{\Pi}^\Phi(y-x) dy.
		\end{align}
		The right-hand side in \eqref{R-LCaputo} can be regarded as a Caputo-D\v{z}rba\v{s}jan (type) derivative. A further relevant fact for \eqref{R-Lxinf} is the equivalence with \eqref{Marchaudleft}. If \eqref{R-LCaputo} holds, then
		\begin{align}
			\label{MarchaudR-L}
			\mathbf{D}_{x-}^\Phi u(x)=\mathcal{D}_{(x,\infty)}^\Phi u(x).
		\end{align}
		Indeed,	using first \eqref{R-LCaputo} and then the second formula in \eqref{tailSymb}, we have
		\begin{align*}
			\mathcal{D}_{(x,\infty)}^\Phi u(x)&=- \int_x^{\infty} u^\prime(y) \overline{\Pi}^\Phi(y-x) dy\\
			&=-\int_0^\infty \frac{d}{dy}u(x+y) \overline{\Pi}^\Phi(y) dy\\
			&=- \int_0^\infty \frac{d}{dy}u(x+y) \Pi^\Phi(y,\infty) dy\\
			&{=}-\int_0^\infty \int_0^z \frac{d}{dy} u(x+y) dy \Pi^\Phi(dz)\\
			&=-\int_0^\infty (u(x+z)-u(x)) \Pi^\Phi(dz)\\
			&=\mathbf{D}_{x-}^\Phi u(x).
		\end{align*}

		Concerning \eqref{R-L0x}, the analogous result for $\mathbf{D}_{x+}^\Phi$ and  $\mathcal{D}_{(0,x)}^\Phi$ can be proved, a clue
		is due to their Laplace transforms.
		
		For the operators \eqref{Marchaudleft} and \eqref{Marchaudright} we introduce the following integration by parts formula.
		\begin{theorem}
			If $u,v \in W_0^{1,1}(0,\infty)$ and $\mathbf{D}_{x-}^\Phi v(x),\mathbf{D}_{x+}^\Phi u(x) \in L^1(0,\infty)$, then
			\begin{align}
				\label{MarchaudAdjoint}
				\int_{0}^\infty u(x) \left(\mathbf{D}_{x-}^\Phi v(x)\right) dx=\int_{0}^\infty \left(\mathbf{D}_{x+}^\Phi u(x)\right)v(x) dx.
			\end{align}
		\end{theorem}
		\begin{proof}
			First of all, from H\"{o}lder's inequality, we observe that
			\begin{align*}
				\vert \vert u \mathbf{D}_{x-}^\Phi v \vert \vert_1 \leq \vert \vert u \vert \vert_\infty \vert \vert \mathbf{D}_{x-}^\Phi v \vert \vert_1 < \infty
			\end{align*}
			because $\mathbf{D}_{x-}^\Phi v \in L^1(0,\infty)$ and $W_0^{1,1}(0,\infty)$ embeds into $L^\infty(0,\infty)$ (see \cite[section 11.2]{leoni2017first}). Similarly for $\mathbf{D}_{x+}$. By Fubini's theorem, from the definition \eqref{Marchaudleft},
			\begin{align*}
				&\int_{0}^\infty u(x) \left(\mathbf{D}_{x-}^\Phi v(x)\right) dx\\ &=\int_{0}^\infty \int_0^\infty u(x)[v(x) - v(x+y)] \Pi^\Phi(dy) dx\\
				&=\int_{0}^\infty \int_0^\infty [u(x)-u(x-y) + u(x-y)][v(x) - v(x+y)] \Pi^\Phi(dy) dx.
			\end{align*}
			By using \eqref{Marchaudright},
			\begin{align*}
				\int_{0}^\infty \int_0^\infty [u(x)-u(x-y)] v(x) \Pi^\Phi(dy) dx= \int_{0}^\infty \left(\mathbf{D}_{x+}^\Phi u(x)\right)v(x) dx,
			\end{align*}
			and
			\begin{align*}
				&\int_{0}^\infty u(x) \left(\mathbf{D}_{x-}^\Phi v(x)\right) dx=\int_{0}^\infty \left(\mathbf{D}_{x+}^\Phi u(x)\right)v(x) dx \ + F(u,v)
			\end{align*}
			where
			\begin{align*}
				F(u,v)&:= \int_{0}^\infty \int_0^\infty \left( u(x-y)[v(x) - v(x+y)] - [u(x)-u(x-y)] v(x+y) \right )\Pi^\Phi(dy) dx\\
				&=\int_{0}^\infty \int_y^\infty u(x-y) v(x) dx \Pi^\Phi(dy)  - \int_{0}^\infty \int_0^\infty u(x) v(x+y) \Pi^\Phi(dy) dx\\
				&=\int_{0}^\infty \int_{0}^\infty u(x) v(x+y)  dx \Pi^\Phi(dy)- \int_{0}^\infty \int_0^\infty u(x) v(x+y) \Pi^\Phi(dy) dx\\
				&= \int_{0}^\infty \int_0^\infty u(x) v(x+y) \Pi^\Phi(dy) dx- \int_{0}^\infty \int_0^\infty u(x) v(x+y) \Pi^\Phi(dy) dx=0,
			\end{align*}
			hence \eqref{MarchaudAdjoint} holds.
		\end{proof}
		\begin{remark}
			The same result for left and right Marchaud derivatives $\mathbf{D}_{x-}^\alpha$ and $\mathbf{D}_{x+}^\alpha$ when $x \in \mathbb{R}$ is presented in \cite[(6.27)]{samko1993fractional} and in \cite[Exercise 1.8.2]{kolokoltsov2011markov}.
		\end{remark}
		The following proposition is a simple reformulation of \cite[Formula 6, Section 15]{ito1963brownian}, based on the definitions of subordinators and Bernstein functions.
		\begin{proposition}
			Under the hypothesis of Theorem \ref{fracdyntm}, we have
			\begin{align*} \mathbf{E}_0 \left[ \int _{0}^\infty e^{-\lambda t}f(B_t^\bullet ) \textrm{d}t\right] = \frac{\int _0^\infty R_\lambda ^D f(y) \, \Pi ^\Phi (\textrm{d}y)}{\Phi (\sqrt{\lambda })}, \end{align*}
		\end{proposition}
		
		\begin{proof}
			In the definition of \eqref{process} we have the process $H^\Phi L^\Phi$, which is a right-continuous process with jumps (given the way the subordinator $H^\Phi$ is constructed). Denote by $(t_n)_{n\geq 1}$ the jump times of $H^\Phi$, which are countable by \cite[Theorem 21.3]{sato1999levy}, and by
			\[
			l_n = \Delta H^\Phi_{t_n}=H^\Phi_{t_n}-H^\Phi_{t_n-}
			\]
			the corresponding jump sizes.
			For each jump we define
			\[
			l_n^- := H^\Phi_{t_n-},
			\qquad
			l_n^+ := H^\Phi_{t_n} = H^\Phi_{t_n-} + l_n.
			\]
			Then the half-line can be decomposed into two disjoint parts:
			\begin{align*}
				\mathcal{J}^c &:= \{\, t \geq 0 : H^\Phi_{L^\Phi_t} = t \,\}, \\
				\mathcal{J} &:= \bigcup_{n \geq 1} \mathcal{J}_n
				= \bigcup_{n \geq 1} [\, l_n^-,\, l_n^+ \,).
			\end{align*}
			Here $\mathcal{J}$ is the union of the intervals generated by the jumps of $H^\Phi$. We have to include the point $l_n^-$ because $H^\Phi L^\Phi$ is right-continuous. The process $B^\bullet$ then takes the form
			\begin{align*}
				B^\bullet_t =
				\begin{cases}
					B^+_t &\gamma_t \in \mathcal{J}^c,\\[4pt]
					l_n^+ - \gamma_t + B^+_t &\gamma_t \in \mathcal{J}_n.
				\end{cases}
			\end{align*}
			Indeed, if $\gamma_t \in \mathcal{J}_n$ we have $H^\Phi L^\Phi \gamma = l_n^+$. Note that $\gamma_t \in [l_n^-,l_n^+)$ if and only if $t \in [\gamma_{-}^{-1}(l_n^-), \gamma_{-}^{-1}(l_n^+))$. Since $\gamma_t$ is continuous and its right inverse $\gamma^{-1}$ is right-continuous with jumps, to ensure the inclusion of the point $l_n^-$ we introduce the left-continuous inverse
			\[
			\gamma_{-}^{-1}(t):=\inf \{s \geq 0 : \gamma_s \geq t\}.
			\]
			
			We prefer writing $\gamma_{-}^{-1}(t)$ over $\gamma_{t-}^{-1}$, to avoid confusion with right and left points of $l_n$. Then, for the resolvent of the process, we have
			\begin{align*}
				&\mathbf{E}_{0}  \left[\int_0^\infty e^{-\lambda t} f(B^\bullet_t)\,dt\right]\\
				&= \mathbf{E}_{0}  \left[\int_{\mathcal{J}\cup \mathcal{J}^c} e^{-\lambda t} f(B^\bullet_t)\,dt\right]\\
				&=\mathbf{E}_{0}  \left[\int_{\mathcal{J}} e^{-\lambda t} f(B^\bullet_t)\,dt\right] + \mathbf{E}_{0}  \left[\int_{\mathcal{J}^c} e^{-\lambda t} f(B^\bullet_t)\,dt\right]\\
				&=\sum_{n \geq 1} \mathbf{E}_{0}  \left[\int_{\gamma_{-}^{-1}(l_n^-)}^{\gamma_{-}^{-1}(l_n^+)} e^{-\lambda t} f(B^\bullet_t)\,dt\right] + \mathbf{E}_{0}  \left[\int_{\mathcal{J}^c} e^{-\lambda t} f(B^\bullet_t)\,dt\right],
			\end{align*}
			but we are dealing only with pure-jump subordinators (with drift zero), hence $\mathcal{J}^c$ has zero Lebesgue measure. For the other part, we get
			\begin{align*}
				&\mathbf{E}_{0}  \left[\int_0^\infty e^{-\lambda t} f(B^\bullet_t)\,dt\right]\\
				&=\sum_{n \geq 1} \mathbf{E}_{0}  \left[\int_{\gamma_{-}^{-1}(l_n^-)}^{\gamma_{-}^{-1}(l_n^+)} e^{-\lambda t} f(B^\bullet_t)\,dt\right]\\
				&=\sum_{n \geq 1}  \mathbf{E}_{0} \left[\int_0^{\gamma_{-}^{-1}(l_n^+)-\gamma_{-}^{-1}(l_n^-)} e^{-\lambda (s + \gamma_{-}^{-1}(l_n^-))} f\big(B^\bullet_{s + \gamma_{-}^{-1}(l_n^-)}\big)\,ds\right].
			\end{align*}
			
			Prior to commencing integration, it is necessary to note certain aspects concerning the local time $\gamma$ and its left inverse $\gamma_{-}^{-1}$. Since $\gamma$ is continuous, we have $\gamma_{\gamma^{-1}_{-}(l_n^-)}=l_n^-$. For the left inverse, serving as the upper bound within the integral, we know that (see \cite[Proposition 1.3, Chapter X]{revuz-yor})
			\[
			\gamma^{-1}_{-}(l_n^+) - \gamma_{-}^{-1}(l_n^-)
			= \gamma_{-}^{-1}(l_n) \circ \theta_{\gamma^{-1}_{-}(l_n^-)} ,
			\]
			where \(l_n=l_n^+ - l_n^-\) and \(\theta\) is the shift operator. We also recall the additive functional property for the local time:
			\[
			\gamma_{s+\gamma^{-1}_{-}(l_n^-)}=\gamma_s \circ \theta_{\gamma^{-1}_{-}(l_n^-)} + \gamma_{\gamma^{-1}_{-}(l_n^-)}
			= \gamma_s \circ \theta_{\gamma^{-1}_{-}(l_n^-)} + l_n^- .
			\]
			
			Then, we rewrite
			\begin{align*}
				&\mathbf{E}_{0}  \left[\int_0^\infty e^{-\lambda t} f(B^\bullet_t)\,dt\right]\\
				&=\sum_{n \geq 1}  \mathbf{E}_{0} \left[\int_0^{\gamma_{-}^{-1}(l_n^+)-\gamma_{-}^{-1}(l_n^-)} e^{-\lambda (s + \gamma_{-}^{-1}(l_n^-))} f\big(B^\bullet_{s + \gamma_{-}^{-1}(l_n^-)}\big)\,ds\right]\\
				&=\sum_{n \geq 1} \mathbf{E}_{0}\!\left[ \mathbf{E}_{0}\!\left[ \Big(\int_0^{\gamma_{-}^{-1}(l_n)} e^{-\lambda s } f\big(l_n + B^+_s - \gamma_s\big)\,ds\Big)\circ \theta_{\gamma^{-1}_{-}(l_n^-)} \Big\vert \mathcal{F}^+_{\theta_{\gamma_{-}^{-1}(l_n^-)}}\right]\right],
			\end{align*}
			where \(\mathcal{F}^+\) is the natural filtration of \(B^+\). By the strong Markov property of \(B^+\) at \(\gamma_{-}^{-1}(l_n^-)\) and the fact \(\gamma_{-}^{-1}(l_n^-)= \gamma^{-1}_{l_n^-}\) a.s., we obtain
			\begin{align*}
				&\mathbf{E}_{0}  \left[\int_0^\infty e^{-\lambda t} f(B^\bullet_t)\,dt\right]\\
				&=\sum_{n \geq 1} \mathbf{E}_{0}\!\left[e^{-\lambda \gamma^{-1}_{l_n^-}} \; \mathbf{E}_{0}\!\left[\int_0^{\gamma_{-}^{-1}(l_n)} e^{-\lambda s } f\big(l_n + B^+_s - \gamma_s\big)\,ds \right] \right].
			\end{align*}
			
			The process \(l_n - \gamma_t + B^+\), since \(t \leq \gamma^{-1}(l_n)\), behaves like Brownian motion \(\{B_t, t \leq \tau_0\}\) started at \(l_n\), with \(\tau_0\) the hitting time of \(0\). Hence
			\begin{align*}
				&\mathbf{E}_{0}\!\left[\int_0^\infty e^{-\lambda t} f(B^\bullet_t)\,dt\right]\\
				&= \sum_{n \geq 1} 
				\mathbf{E}_{0}\!\left[
				e^{-\lambda \gamma^{-1}_{t_n^-}}
				\mathbf{E}_{l_n}\!\left[\int_{0}^{\tau_0} e^{-\lambda t} f(B_t)\,dt \right]
				\right]\\
				&= \sum_{n \geq 1} 
				\mathbf{E}_{0}\!\left[
				e^{-\sqrt{\lambda}\, t_n^-}\;
				R_\lambda^D f(l_n)
				\right]\\
				&= \sum_{n \geq 1} 
				\mathbf{E}_{0}\!\left[
				e^{-\sqrt{\lambda}\, H^\Phi_{t_n^-}}\;
				R_\lambda^D f(l_n)
				\right]\\
				&= \mathbf{E}_{0}\!\left[
				\int_{(0,\infty)\times(0,\infty)}
				e^{-\sqrt{\lambda}\, H^\Phi_{t^-}}\;
				R_\lambda^D f(l)\,
				N(dt\times dl)
				\right],
			\end{align*}
			where \((t_n)_{n\ge1}\) are the jump times of the subordinator \(H^\Phi\), \(l_n=\Delta H^\Phi_{t_n}\) the corresponding jump sizes, and \(N(dt\times dl)\) is the Poisson random measure of jumps of \(H^\Phi\).
			
			From \cite[Example II.4.1]{ikeda2014stochastic} we know that \(N(dt \times dl)= dt \, \Pi^\Phi(dl)\), thus
			\begin{align*}
				&\mathbf{E}_{0} \left[\int_{(0,\infty) \times (0,\infty)} e^{-\sqrt{\lambda} H^\Phi_{t^-}} R_\lambda^D f(l)\, N(dt \times dl)\right]\\
				&=\int_0^\infty \int_0^\infty e^{-t \Phi(\sqrt{\lambda})} R_\lambda^D f(l)\, dt \,\Pi^\Phi(dl)\\
				&=\frac{\displaystyle\int_0^\infty R_\lambda^D f(l)\,\Pi^\Phi(dl)}{\Phi(\sqrt{\lambda})}.
			\end{align*}
		\end{proof}

		\section*{Acknowledgments}
		The authors thank the University of Trento, Sapienza University of Rome and the group INdAM-GNAMPA for the support under their Grants.\\
		GP is supported by the Basque Government through the BERC 2022--2025 program and by the Ministry of Science and Innovation: BCAM Severo Ochoa accreditation CEX2021-001142-S / MICIN / AEI / 10.13039/501100011033.\\
		The research has been mostly funded by MUR under the project PRIN 2022 - 2022XZSAFN: Anomalous Phenomena on Regular and Irregular Domains: Approximating Complexity for the Applied Sciences - CUP B53D23009540006 - M4.C2.1.1 - Web Site: \url{https://www.sbai.uniroma1.it/~mirko.dovidio/prinSite/index.html}.

		\medskip

\begin{thebibliography}{99}
			
			\bibitem{yor-walsh}
			Barlow, M., Pitman, J., Yor, M.:
			On Walsh's Brownian motions. In: S\'eminaire de Probabilit\'es, XXIII, Lecture Notes in Math., vol. 1372, pp. 275--293. Springer, Berlin (1989).
			
			\bibitem{berry2025stationary}
			Berry, J., Camilli, F.:
			Stationary Mean Field Games on networks with sticky transition conditions. ESAIM: Control Optim. Calc. Var. \textbf{31}, 18 (2025).
			
			\bibitem{berry2024sticky}
			Berry, J., Colantoni, F.:
			Sticky diffusions on star graphs: characterization and It\^o formula. arXiv:2411.05441 (2024).
			
			\bibitem{bertoin1996levy}
			Bertoin, J.:
			L\'evy processes. Cambridge University Press, Cambridge (1996).
			
			\bibitem{bertoin1999subordinators}
			Bertoin, J.:
			Subordinators: examples and applications. In: Lectures on probability theory and statistics, pp. 1--91. Springer (1999).
			
			\bibitem{blumenthal}
			Blumenthal, R. M., Getoor, R. K.:
			Markov processes and potential theory. Academic Press, New York--London (1968).
			
			\bibitem{Bogdan}
			Bogdan, K., Kunze, M.:
			Stable processes with reflections. Stoch. Process. Appl. \textbf{187}, 104654 (2025).
			
			\bibitem{BonDov}
			Bonaccorsi, S., D'Ovidio, M.:
			Sticky Brownian motions on star graphs.
			Fract. Calc. Appl. Anal. \textbf{27}, 2859--2891 (2024).
			
			\bibitem{handbook}
			Borodin, A. N., Salminen, P.:
			Handbook of Brownian motion—facts and formulae, 2nd edn. Birkh\"auser, Basel (2002).
			
			\bibitem{BKLT-financeSticky}
			Bouchaud, J., Kr\"uger, P., Landier, A., Thesmar, D.:
			Sticky Expectations and the Profitability Anomaly. J. Finance \textbf{74}, 639--674 (2018).
			
			\bibitem{capitanelli2019fractional}
			Capitanelli, R., D'Ovidio, M.:
			Fractional equations via convergence of forms. Fract. Calc. Appl. Anal. \textbf{22}(4), 844--870 (2019).
			
			\bibitem{capitanelli2020delayed}
			Capitanelli, R., D'Ovidio, M.:
			Delayed and rushed motions through time change. ALEA \textbf{17}, 183--204 (2020).
			
			\bibitem{caputoBook}
			Caputo, M.:
			Elasticit\`a e Dissipazione. Zanichelli, Bologna (1969).
			
			\bibitem{CapMai71}
			Caputo, M., Mainardi, F.:
			Linear models of dissipation in anelastic solids. La Rivista del Nuovo Cimento \textbf{1}, 161--198 (1971).
			
			\bibitem{CapMai71b}
			Caputo, M., Mainardi, F.:
			A new dissipation model based on memory mechanism. PAGEOPH \textbf{91}, 134--147 (1971).
			
			\bibitem{chen2017time}
			Chen, Z.-Q.:
			Time fractional equations and probabilistic representation. Chaos Solitons Fractals \textbf{102}, 168--174 (2017).
			
			\bibitem{chung}
			Chung, K. L.:
			Excursions in Brownian motion. Ark. Mat. \textbf{14}(2), 155--177 (1976).
			
			\bibitem{colantoni2025non}
			Colantoni, F.:
			Non-local skew and non-local skew sticky Brownian motions.
			J. Evol. Equ. \textbf{25}(2), 1--26 (2025).
			
			\bibitem{colantoni2021inverse}
			Colantoni, F., D'Ovidio, M.:
			On the inverse gamma subordinator. Stochastic Anal. Appl. \textbf{41}(5), 999--1024 (2023).
			
			\bibitem{ColDovPag}
			Colantoni, F., D'Ovidio, M., Pagnini, G.:
			Time reversal of Brownian motion with Poissonian resetting. arXiv:2505.15639 (2025). To appear in J. Stat. Phys.
			
			\bibitem{ColDovTav}
			Colantoni, F., D'Ovidio, M., Tavani, F.:
			Earthquake modelling via Brownian motions on networks. arXiv:2509.21851 (2025).
			
			\bibitem{d2022fractional}
			D'Ovidio, M.:
			Fractional boundary value problems. Fract. Calc. Appl. Anal. 	\textbf{25}(1), 29--59 (2022).
			
			\bibitem{d2022fractionalsticky}
			D'Ovidio, M.:
			Fractional boundary value problems and elastic sticky Brownian motions. Fract. Calc. Appl. Anal. \textbf{27}, 2162--2202 (2024).
			
			\bibitem{d2024fractionalsticky}
			D'Ovidio, M.:
			Fractional boundary value problems and elastic sticky Brownian motions, II: The bounded domain. arXiv:2205.04162 (2024).
			
			\bibitem{Dzh66}
			D\v{z}rba\v{s}jan, M. M.:
			Integral transforms and representations of functions in the complex plane. Nauka, Moscow (1966) (in Russian).
			
			\bibitem{DzhNers68}
			D\v{z}rba\v{s}jan, M. M., Nersessian, A. B.:
			Fractional derivatives and the Cauchy problem for differential equations of fractional order. Izv. Akad. Nauk Armjan. SSR. Ser. Mat. \textbf{3} (1968), No. 1, 1--29 (in Russian).
			
			\bibitem{evans2010partial}
			Evans, L. C.:
			Partial differential equations. American Mathematical Society, Providence (2010).
			
			\bibitem{resetting1}
			Evans, M. R., Majumdar, S. N.:
			Diffusion with stochastic resetting. Phys. Rev. Lett. \textbf{106}(16), 160601 (2011).
			
			\bibitem{resetting2}
			Evans, M. R., Majumdar, S. N., Schehr, G.:
			Stochastic resetting and applications. J. Phys. A: Math. Theor. \textbf{53}, 193001 (2020).
			
			\bibitem{feller1952parabolic}
			Feller, W.:
			The parabolic differential equations and the associated semi-groups of transformations. Ann. Math. \textbf{55}(3), 468--519 (1952).
			
			\bibitem{ferrari2018weyl}
			Ferrari, F.:
			Weyl and Marchaud derivatives: A forgotten history. Mathematics \textbf{6}(1), 6 (2018).
			
			\bibitem{GNSS-finance-Sticky}
			Gerali, A., Neri, S., Sessa, L., Signoretti, F. M.:
			Credit and Banking in a DSGE Model of the Euro Area. J. Money Credit Bank. \textbf{42}, Suppl. 1, 107--141 (2010).
			
			\bibitem{HSU}
			Hsu, P.:
			On excursions of reflecting Brownian motion. Trans. Amer. Math. Soc. \textbf{296}(1), 239--264 (1986).
			
			\bibitem{ikeda2014stochastic}
			Ikeda, N., Watanabe, S.:
			Stochastic differential equations and diffusion processes. Elsevier (2014).
			
			\bibitem{ito1963brownian}
			It\^{o}, K., McKean, H. P., Jr.:
			Brownian motions on a half line. Illinois J. Math. \textbf{7}(2), 181--231 (1963).
			
			\bibitem{ito1996diffusion}
			It\^{o}, K., McKean, H. P., Jr.:
			Diffusion processes and their sample paths. Springer-Verlag, Berlin--New York (1974). Second printing, corrected.
			
			\bibitem{sato1999levy}
			Sato, K.-I.:
			L\'evy processes and infinitely divisible distributions. Cambridge University Press, Cambridge (1999).
			
			\bibitem{kilbas2006theory}
			Kilbas, A. A., Srivastava, H. M., Trujillo, J. J.:
			Theory and applications of fractional differential equations, North-Holland Math. Stud., vol. 204. Elsevier Science B.V., Amsterdam (2006).
			
			\bibitem{kochubei2011general}
			Kochubei, A. N.:
			General fractional calculus, evolution equations, and renewal processes. Integral Equ. Oper. Theory \textbf{71}(4), 583--600 (2011).
			
			\bibitem{russi3}
			Kostrykin, V., Potthoff, J., Schrader, R.:
			Brownian motions on metric graphs. J. Math. Phys. \textbf{53}(9), 095206 (2012).
			
			\bibitem{russi1}
			Kostrykin, V., Potthoff, J., Schrader, R.:
			Construction of the paths of Brownian motions on star graphs I. Commun. Stoch. Anal. \textbf{6}(2), 223--245 (2012).
			
			\bibitem{russi2}
			Kostrykin, V., Potthoff, J., Schrader, R.:
			Construction of the paths of Brownian motions on star graphs II. Commun. Stoch. Anal. \textbf{6}(2), 247--261 (2012).
			
			\bibitem{kolokoltsov2011markov}
			Kolokoltsov, V. N.:
			Markov processes, semigroups and generators. De Gruyter Stud. Math., vol. 38. Walter de Gruyter \& Co., Berlin (2011).
			
			\bibitem{leoni2017first}
			Leoni, G.:
			A first course in Sobolev spaces. American Mathematical Society, Providence (2017).
			
			\bibitem{MaLuPa2001}
			Mainardi, F., Luchko, Y., Pagnini, G.:
			The fundamental solution of the space-time fractional diffusion equation. Fract. Calc. Appl. Anal. \textbf{4}(2), 153--192 (2001).
			
			\bibitem{Mais1983}
			Maisonneuve, B.:
			Ensembles r\'eg\'en\'eratifs de la droite. Z. Wahrscheinlichkeitstheorie verw. Gebiete \textbf{63}, 501--510 (1983).
			
			\bibitem{MNV09}
			Meerschaert, M. M., Nane, E., Vellaisamy, P.:
			Fractional Cauchy problems on bounded domains. Ann. Probab. \textbf{37}(3), 979--1007 (2009).
			
			\bibitem{MeSt14}
			Meerschaert, M. M., Straka, P.:
			Semi-Markov approach to continuous time random walk limit processes. Ann. Probab. \textbf{42}(4), 1699--1723 (2014).
			
			\bibitem{meerschaert2008triangular}
			Meerschaert, M. M., Scheffler, H.-P.:
			Triangular array limits for continuous time random walks. Stoch. Process. Appl. \textbf{118}(9), 1606--1633 (2008).
			
			\bibitem{meerschaert2019stochastic}
			Meerschaert, M. M., Sikorskii, A.:
			Stochastic models for fractional calculus. De Gruyter (2019).
			
			\bibitem{bookMorPer}
			M\"orters, P., Peres, Y.:
			Brownian Motion. Cambridge University Press (2012).
			
			\bibitem{mugnolo-libro}
			Mugnolo, D.:
			Semigroup methods for evolution equations on networks. Springer, Cham (2014).
			
			\bibitem{mugnolo2019actually}
			Mugnolo, D.:
			What is actually a metric graph? arXiv:1912.07549 (2019).
			
			\bibitem{phillips1952generation}
			Phillips, R. S.:
			On the generation of semigroups of linear operators. Pacific J. Math. \textbf{2}(3), 343--369 (1952).
			
			\bibitem{revuz-yor}
			Revuz, D., Yor, M.:
			Continuous martingales and Brownian motion, 3rd edn. Springer, Berlin (1999).
			
			\bibitem{samko1993fractional}
			Samko, S. G., Kilbas, A. A., Marichev, O. I.:
			Fractional integrals and derivatives. Gordon and Breach Science Publishers, Yverdon (1993).
			
			\bibitem{schilling2012bernstein}
			Schilling, R. L., Song, R., Vondracek, Z.:
			Bernstein functions. De Gruyter (2012).
			
			\bibitem{toaldo2015convolution}
			Toaldo, B.:
			Convolution-type derivatives, hitting-times of subordinators and time-changed $C_0$-semigroups. Potential Anal. 	\textbf{42}(1), 115--140 (2015).
			
			\bibitem{walsh1978diffusion}
			Walsh, J. B.:
			A diffusion with a discontinuous local time. Ast\'erisque \textbf{52}(53), 37--45 (1978).
			
			\bibitem{werner2016brownian}
			Werner, F.:
			Brownian motions on metric graphs. PhD thesis, University of Mannheim (2016).
			
			\bibitem{widder1944positive}
			Widder, D. V.:
			Positive temperatures on an infinite rod. Trans. Amer. Math. Soc. \textbf{55}, 85--95 (1944).
			
		\end{thebibliography}
		
\end{document}